\documentclass{article}

\usepackage[utf8]{inputenc}
\usepackage{geometry}
\usepackage{amsmath}
\usepackage{amssymb}
\usepackage{amsthm}
\usepackage{graphicx}
\usepackage{caption}
\usepackage{subcaption}
\usepackage{listings}
\usepackage{color}
\usepackage{float}
\usepackage{siunitx}
\usepackage{placeins}
\usepackage{booktabs}
\usepackage{enumitem}
\usepackage{paralist}
\usepackage{hyperref}

\usepackage[backend = biber,
    style = alphabetic,
    giveninits = true,
    natbib = true,
    hyperref = true,
    maxbibnames = 10,
    ]{biblatex}

\usepackage{algorithm}
\usepackage[noend]{algpseudocode}

\makeatletter
\let\@fnsymbol\@arabic
\makeatother

\definecolor{Mat1}{rgb}{     0,0.4470,0.7410}
\definecolor{Mat2}{rgb}{0.8500,0.3250,0.0980}
\definecolor{Mat3}{rgb}{0.9290,0.6940,0.1250}
\definecolor{Mat4}{rgb}{0.4940,0.1840,0.5560}
\definecolor{Mat5}{rgb}{0.4660,0.6740,0.1880}
\definecolor{Mat6}{rgb}{0.3010,0.7450,0.9330}
\definecolor{Mat7}{rgb}{0.6350,0.0780,0.1840}

\newcommand{\tr}{ {\mathsf T} }
\newcommand{\ct}{ {\mathsf H} }
\newcommand{\ima}{\operatorname{i}}
\newcommand{\pp}{\operatorname{p}}

\newcommand{\cc}{c.c.}
\newcommand{\env}{\operatorname{env}}

\newcommand{\para}{q}
\newcommand{\vpara}{\mathbf{q}}
\newcommand{\paraspace}{\boldsymbol{Q}} % parameter space
\newcommand{\Fo}{\mathbf{F}}
\newcommand{\Nsample}{M}
\newcommand{\CC}{\mathbb{C}}
\newcommand{\NN}{\mathbb{N}}
\newcommand{\ve}{\mathbf{e}}
\newcommand{\vf}{\mathbf{f}}
\newcommand{\vn}{\mathbf{n}}
\newcommand{\vo}{\mathbf{0}}
\newcommand{\vu}{\mathbf{u}}

\newcommand{\vx}{\mathbf{x}}
\newcommand{\vy}{\mathbf{y}}
\newcommand{\mA}{\mathbf{A}}
\newcommand{\mC}{\mathbf{C}}
\newcommand{\mD}{\mathbf{D}}

\newcommand{\mI}{\mathbf{I}}

\newcommand{\mS}{\mathbf{S}}
\newcommand{\mV}{\mathbf{V}}
\newcommand{\mX}{\mathbf{X}}

\newcommand{\vsigma}{\pmb{\sigma}}
\newcommand{\vnabla}{\pmb{\nabla}}

\newcommand{\ff}{\boldsymbol{f}}
\newcommand{\fu}{\boldsymbol{u}}

\newcommand{\dOmega}{\omega_\Delta}
\newcommand{\ecf}{f_c}      % excitation center frequenz 
\newcommand{\ect}{t_c}      % excitation center time shift 
\newcommand{\EWM}{\zeta}    % EWM complex part
\newcommand{\cg}{c_g}       % group 
\newcommand{\Ldefect}{L_1}  % length difference
\newcommand{\Tdiff}{T_{\Delta}}  % time difference
\newcommand{\Tex}{T_{\text{ex}}}  % time difference
\newcommand{\cgP}{c_g^{\text{ex}}}  % group velocity exicitation
\newcommand{\cgN}{c_g^{\text{re}}}  % group velocity reflection
\definecolor{darkgreen}{rgb}{0,0.6,0}

\title{Defect  reconstruction  in  a  2D  semi-analytical  waveguide  model via derivative-based optimization}

\date{\today}

\author{J. Bulling\thanks{Bundesanstalt für Materialforschung und -pruefung, Unter den Eichen 87, 12205 Berlin, Germany \newline Jannis.Bullling@BAM.de, Jens.Prager@BAM.de} \and B. Jurgelucks\thanks{Department of Mathematics, Humboldt-Universität zu Berlin, Unter den Linden 6, 10099 Berlin, Germany\newline Benjamin.Jurgelucks@math.hu-berlin.de, Andrea.Walther@math.hu-berlin.de}, \and J. Prager$^1$, \and A. Walther$^2$}

\begin{document}
    \maketitle

	This paper considers the reconstruction of a defect in a two-dimensional waveguide during non-destructive ultrasonic inspection using a derivative-based optimization approach. The propagation of the mechanical waves is simulated by the Scaled Boundary Finite Element Method (SBFEM) that builds on a semi-analytical approach. The simulated data is then fitted to a given set of data describing the reflection of a defect to be reconstructed. For this purpose, we apply an iteratively regularized Gauss-Newton method in combination with algorithmic differentiation to provide the required derivative information accurately and efficiently. We present numerical results for three different kinds of defects, namely a crack, a delamination, and a corrosion. These examples show that the parameterization of the defect can be reconstructed efficiently and robustly in the presence of noise.
	\section{Introduction}

Ultrasonic guided waves, especially Lamb-waves, are used in Non-Destructive Testing (NDT) and Structural Health Monitoring (SHM) to identify defects. Here the long testing range that can be covered using guided waves is one of the main advantages over conventional scanning techniques. However, this advantage also yields new challenges. The multi-modal nature of guided waves with different wave velocities and the large distance between the defect and the sensor makes the localisation and characterisation of defects more difficult. On the other side, Lamb-waves provide the flexibility to choose the most appropriate mode, which is the most sensitive to a specific type of defect.

\begin{figure}[t]
    \begin{center}
    	\includegraphics[width=.4\textwidth]{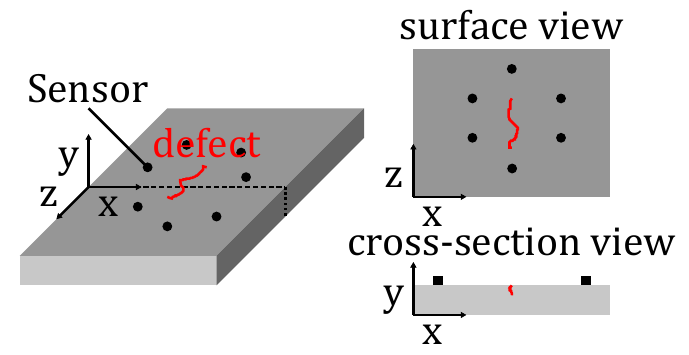}
		\caption{Different views for a 2d-simplification of the 3d-problem}
    	\label{fig:View}
    \end{center}
\end{figure}

Most models for crack identification assume that the data containing the waves scattered by the defect are available at certain sensor positions. Depending on the type of sensor, this data can be an electrical signal that represents a displacement or velocity field. We will refer to this data as measurement data, even though for many publications including this one, it is generated artificially by appropriate simulations. There are two main types of methods to detect a defect in a waveguide. One approach is based on imaging to localise and characterise the defects. Another class of algorithms fits a damage model directly with the measured data. The resulting damage model can then be further analyzed, e.g., using fracture mechanics techniques. This also supports other goals like NDT and SHM investigations and may answer the question of whether the damage is severe leading to a failure of the structure.

Most imaging algorithms assume a sensor network around the defect and consider a surface view of the waveguide as sketched in Figure~\ref{fig:View}. Often these approaches are called guided wave tomography or migration in the seismic community. The measurement data is then propagated backwards from the sensors on the surface. As a possible method of backward propagation, each pixel of the surface is analyzed by the direct sound path between this pixel and the transmitter/receiver pairs \cite{ihn2008pitch}. This approach relies on the assumption that there is only one dominant wave mode to use its group velocity for the backward propagation. Other methods do not need a dominating mode, for example, if a time-reversal approach is utilized as backward propagation method~\cite{mori2019damage}. There are many competing approaches for the backward propagation method. Since a full review is out of the scope of this introduction, we refer the reader to \cite{zhao2011ultrasonic}~and~\cite{neubeck2020matrix} for a comparison of the various methods. If the backward propagation method is defined by the transmitter/receiver-path, it depends strongly on the number of sensors and their position~\cite{zhao2011ultrasonic}. In other cases, the resolution of these approaches is directly linked to the physical properties of the propagating waves, for example, the wavelength~\cite{mori2019damage}. Most tomography algorithms show only a projection on the surface of the waveguide and do not show any depth of the defect.

As an alternative to the above imaging algorithms, the localization and characterization of defects can be viewed as an inverse problem. The algorithm tries to find the set of parameters of a damage model that yield the smallest distance of the simulated data to the measurement data. These approaches consist of three main components. First, a forward model that generates simulated data representing the damage for the current set of parameters. Second, one uses an objective function that compares the measurement and simulation. Third, an optimization algorithm updates the parameter values  to minimize the objective function, i.e., obtain a better fit of simulation and measurement. With regard to the simulation method, damage model, target function and optimization, there are different approaches that are briefly discussed in the following paragraphs.

For a discrete damage model, the extended finite element method (XFEM) is often used, because of its possibility to define mesh independent cracks~\cite{rabinovich2007xfem, rabinovich2009crack, jung2014modeling, agathos2018multiple, livani2018identification}. Rabinovich et al.~\cite{rabinovich2009crack} used the XFEM in the frequency domain to reconstruct a single straight crack. The objective function is based on the least square error between measurement and simulation data. This objective function oscillates and has several local minima. A genetic algorithm optimizes the objective function to overcome the difficulty of getting trapped in a local minimum. The same authors improved the results by switching to the time-domain~\cite{rabinovich2009crack} and changing the objective function. The objective function in the time-domain is based on the "arrival time" of the reflection, which is often used in classic defect detection in an ultrasonic test. Using the arrival time, the objective function oscillates less, but there are still local minima. Livani et al. \cite{livani2018identification} investigated multiple defects with particle swarm optimization. Finally, three dimensional problems with multiple cracks are solved by Agathos et al.~\cite{agathos2018multiple}. They used a two-step GA approach for elliptical cracks. Jung et al.~\cite{jung2014modeling} consider curved cracks and inclusions together with a gradient descent as the optimization method. Multiple optimization runs are performed with different starting points representing different cracks to circumvent local minima. 
% thematisch sortiert nach genetischen Algorithmen und Ableitungsbasiert

A time-reversal approach can also be used to optimize defect parameters. If a sensor signal is reversed in time and transmitted back in the domain, it will be refocused at the source. Amitt et al.~\cite{amitt2014time} defined a objective function based on these refocusing of the wave in a specific part of a membrane\footnote{The word membrane should indicate that the Helmholtz equation is solved instead of the elastic wave equation.}. Unfortunately, this objective function is oscillating~\cite{amitt2014time}. A full scan of the parameter space that defined the crack and its position was performed for optimization. \cite{lopatin2017computational} considers an experimental validation of this approach.

Seidle et al.~\cite{seidl2016iterative} used another time-reversal approach for a membrane. Here, an integral damage model is used which lowers the stiffness for certain elements. The connection between XFEM and a stiffness-reduction can be found in~\cite{broumand2021inverse}. This approach, also known as Full Wavefield Inversion, allows many possible defect geometries and is not limiting their number, but the resolution of the defect depends on the computational grid. A highly refined computational grid can lead to very time-consuming calculations, especially for three-dimensional problems. A theoretical and experimental investigation of a similar approach can be found, e.g., in Rao et al.~\cite{rao2016guided,rao2017investigation,rao2020multi}, where an acoustic inversion is used for elastic wave propagation to lower the computational demand. A change in both stiffness and density is considered in~\cite{rao2020multi}.

The approaches mentioned above exploit the surface of the waveguide, whereas also a cross-section view may serve as alternative description. A plane strain simplification can be used to derive a two-dimensional problem - see Figure~\ref{fig:View}. Here, it is important to stress that these alternative dimensions of the cross-sectional view enable the analysis of the depth of the defect and the mode conversion due to defects with a certain depth. Wu et al.~\cite{wu2002inverse} analyzed composite plates with the strip element method. The model allows investigating horizontal and vertical cracks in the cross-section. Semi-analytical methods such as the strip element method have the advantage that the forward calculation is efficient. This simulation approach is combined with a genetic algorithm to  perform the optimization.

Gravenkamp proposed another semi-analytical method, the Scaled Boundary Finite Element Method (SBFEM), for simulating plate-like structures in the context of NDT and SHM~\cite{gravenkamp2012simulation}. The SBFEM was first developed by Song and Wolf to compute bounded and unbounded domains with the same approach~\cite{song2000scaled,wolf2000scaled}. The domain is divided into super elements. Quasi-polar coordinates, called scaled boundary coordinates, construct this first type of super element. For problems in the frequency domain, the super element has only degrees of freedom on the boundary, which leads to a reduction of one dimension in the computational cost. Later, Gravenkamp et al. used another coordinate transformation to construct super elements with a constant cross-section for bounded and unbounded domains in two-dimensions~\cite{gravenkamp2012numerical,gravenkamp2015simulation}. These super elements are very efficient for simulating waveguides. The approach was extended to curved waveguides, three-dimensional-waveguides and waveguides with fluid/structure interaction~\cite{krome2017semi,krome2017prismatic,wasmer2018fluid}. In addition to the simulation of unbounded domains, a discrete crack model is an essential development in SBFEM~\cite{song2018review}.
Quite recently, the SBFEM was used in a shape optimization scheme for a horn speaker and meta-materials~\cite{khajah2021shape} illustrating the flexibility of this formulation in the meshing process and and with respect to infinite boundary conditions.

In this contribution, we propose a general method for reconstructing a single discrete defect model inside a cross-section model of a waveguide. The SBFEM is used as a efficient forward model for the simulation.
The defect  reconstruction is formulated as an optimization problem and a derivative-based optimization algorithm is used to solve it in a systematic and efficient way. The technique of Algorithmic Differentiation (AD, \cite{GrWa08,Nau12}) is used to provide the required gradient information exactly and with low computational costs. We analyse the resulting approach with respect to properties of the objective function, the reconstruction quality and the robustness against noise.

% The problem of defect reconstruction for waveguides is also related to inverse scattering problem~\cite{kress1996inverse,colton2000recent,guzina2004topological}, but the inverse scattering problem considers infinite full and half-spaces. To the best of the author's knowledge, the theory of the inverse scattering problem was not extended to waveguides.

The remainder of this paper is structured in the following way: The forward model and the physical assumptions are provided in Section~\ref{sec:ForwardModel}. Section~\ref{sec:Optimization} summarizes the optimization scheme. This includes a discussion about the objective function and an estimate for the initial value. Three numerical examples are presented in section~\ref{sec:NumericalExperiments} showing the versatility for different damage geometries. Conclusions are given in section~\ref{sec:Conclusion}. The appendix contains additional information on the implementation.

	\section{SBFEM forward model}\label{sec:ForwardModel}

This section introduces the forward model to simulate the guided waves. Throughout, we consider only the two-dimensional case. The extension to three dimensions is subject of our current research. As mentioned in the introduction, the damage is incorporated directly into the waveguide model such that specific parameters model the damage.
The forward model incorporates all important aspects of an ultrasound test on a waveguide. 
However, some aspects can only be designed for a specific experiment. 
Hence, certain simplifying assumptions have to be made to achieve greater general validity.

The first of these assumptions is that the defects are in a specific area within the waveguide modeled as an infinite domain. It is assumed that the vertical edges are far enough away from this area of interest so that the reflections from the vertical edges of the waveguide do not interfere with the reflections from the damage.

The second assumption is that the excitation by a sensor can be modeled by a traction force $\check{\ff}(\vx,t)$ on one part of the boundary.
This traction is mapped into the frequency domain by the Laplace-Transformation.
In contrast to the Fourier-Transformation, the imaginary part of the Laplace-Transformation should weaken the resonance and wrap-around effects of the numerical model. 
In the discrete-time version, this procedure, known as Exponential Window Method~(EWM) \cite{kausel2017advanced}, is given by
\begin{align}
	\check{\ff}(\vx,t) = \ff(\vx) \check{\tau}(t) = \ff(\vx) \sum_{n=0}^N (\tau_{\omega_n} \exp(-\ima \omega_n t) + \cc),
\end{align}
where $\cc$ denotes the complex conjugated of the previous term.
The complex angular  frequency $\omega_n = (n\dOmega - \ima\EWM)\in\CC$ is defined by a uniform frequency step $\dOmega$ and a small parameter $\EWM$.
Note that only the real part of the complex angular frequencies changes.

Algorithmically, the EWM is implemented in the following three steps:
\begin{compactenum}[(i)]
    \item Multiply $\check{\tau}(t)$ with the window function $\exp(-\EWM t)$ with $\EWM = 0.5 \dOmega$, where the frequency step $\dOmega$ results  from the time vector,
    use the Fast-Fourier-Transform (FFT) on the product, and search for the relevant  frequencies.
    \item Approximate the displacement spectrum for all relevant complex frequencies $\omega_n$ (see Equation~\eqref{eq:PDE},~\eqref{eq:SLE}).
    \item For evaluation in the time-domain, transform the displacement by the inverse FFT and afterwards multiply the inverse window function $\exp(\EWM t)$.
\end{compactenum}

\begin{figure}[ht]
    \begin{center}
        \begin{subfigure}{.4\textwidth}
            \centering            \includegraphics[width=1\textwidth]{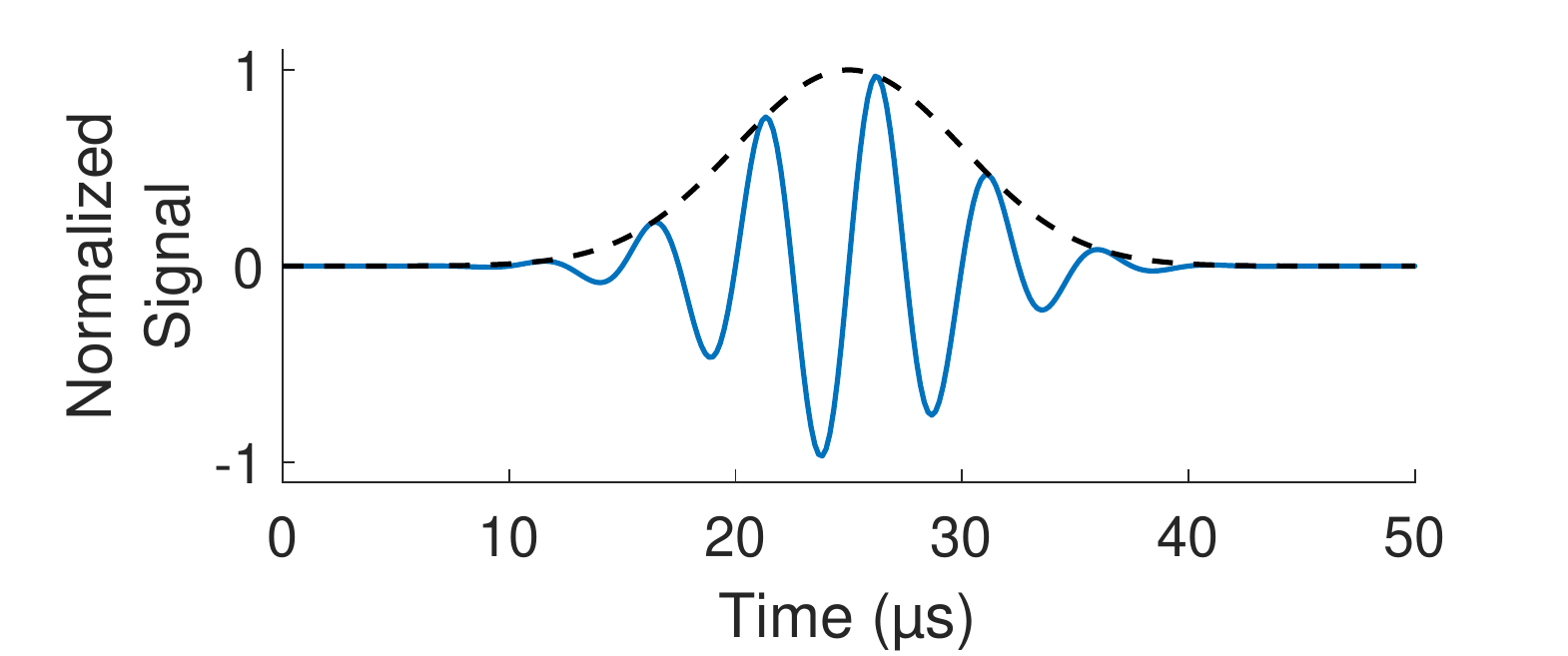}
            \caption{Excitation in the time-domain.}
            \label{fig:excitationA}
        \end{subfigure}
        \begin{subfigure}{.4\textwidth}
            \centering
            \includegraphics[width=1\textwidth]{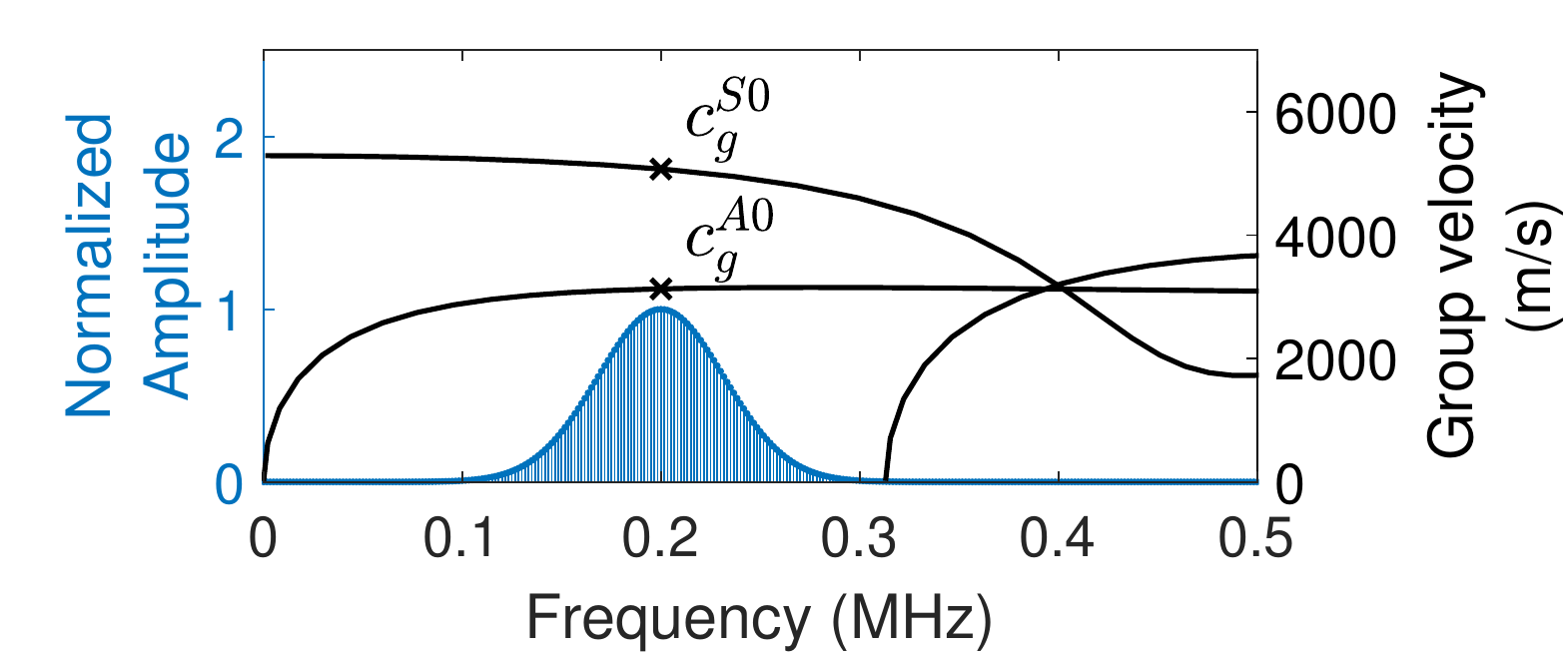}
            \caption{Excitation in the frequency-domain.}
            \label{fig:excitationB}
        \end{subfigure}
        \caption{Excitation.}
        \label{fig:excitation}
    \end{center}
\end{figure}

For all examples, the time dependent part $\check{\tau}(t)$ of the excitation is a sin-modulated Gaussian pulse, i.e.,
\begin{equation}
	\check{\tau}(t) = \sin(2\pi \ecf t)\exp\left(- \frac{(t - \ect)^2}{2\ecf^{-2}}\right)
\end{equation}
with the center frequency $\ecf = \SI{200}{\kilo\hertz}$ and a time shift of $\ect = 5\ecf^{-1}$. Figure~\ref{fig:excitationA} shows the pulse in the time-domain, while \ref{fig:excitationB} illustrates the spectrum. Additionally, the dispersion curves of the waveguide are presented~\cite{gravenkamp2012numerical}. This figure shows that the relevant spectrum of the excitation lies in the area, where only the two fundamental modes are propagating.

For each discrete angular frequency step, the SBFEM is used to approximate the displacement in the equations of linear elastodynamics
\begin{equation}
    \begin{aligned}
    	\vnabla \cdot \vsigma\left(\fu(\vx)\right) + \omega_n^2 \rho  \fu(\vx) & = \vo & \vx & \in \Omega \\
    	\vn^\tr \vsigma\left(\fu(\vx)\right) & = \tau_{\omega_n} \ff(\vx) & \vx & \in \Gamma \label{eq:PDE}
    \end{aligned}
\end{equation}
with the displacement $\vu$, the linear stress $\vsigma$, the density $\rho$ and $\vn$ is the outer normal vector.

\begin{figure}[ht]
    \begin{center}
    \includegraphics[width=1\textwidth]{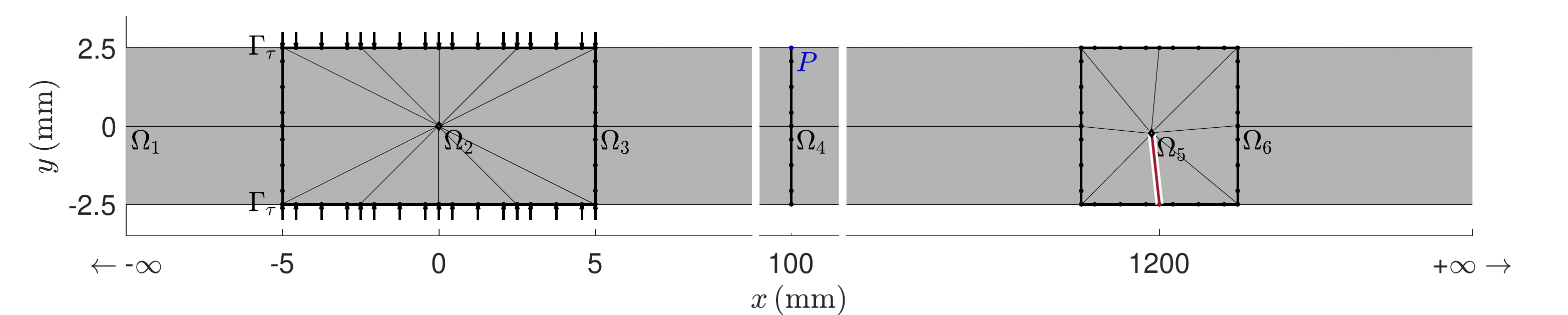}
		\caption{Domain for the \textit{waveguide with a crack}.}
		\label{fig:Crack}
    \end{center}
\end{figure}

As for many other numerical methods, when using the SBFEM, the PDE in Equation~\eqref{eq:PDE} is converted into systems of linear equations
\begin{align}
	\mS_{\omega_n} \vu_{\omega_n} = \vf_{\omega_n} \label{eq:SLE}
\end{align}
with the global, dynamic stiffness matrix $\mS_{\omega_n}$ that depends on the complex angular frequency, the nodal displacement vector $\vu_{\omega_n}$ and the nodal traction vector $\vf_{\omega_n}$.
To compute the global dynamic stiffness matrix $\mS_{\omega_n}$, the domain has to be sub-divided into super elements, i.e., $\Omega = \bigcup_{k = 1}^K \Omega_k$.
Figure~\ref{fig:Crack} shows such a subdivision into super elements.
For each super element $\Omega_k$, parts of the boundary are meshed with finite elements that are used to compute a local stiffness matrix $\mS_{\omega_n}^k$.
These local stiffness matrices are then assembled into the global matrix, while the global nodal traction vector is assembled directly from the finite elements.
The finite element approximation uses spectral shape functions based on the Gauss-Lobatto points.
Figure~\Ref{fig:Crack} depicts elements with shape functions of degree $\pp = 4$ with 5 nodes (marked with $\bullet$).
The corner points of the finite elements are dependent on a set of parameters $\para_i$ that is later modified by the optimization process. 

The prismatic super elements are used to approximate the extended parts of the waveguide.
Figure~\ref{fig:Crack} shows examples with the super elements $\Omega_1$, $\Omega_3$, $\Omega_4$ and $\Omega_6$.
These super elements require a constant cross-section.
The cross-section parallel to the y-axis is then meshed with finite elements and scaled in the $x$-direction.
If the super element is finite, as $\Omega_3$ in Figure~\ref{fig:Crack}, the finite element meshes must coincide on both sides.
Note that the length in $x$-direction does not influence the computational cost because there are only nodes at the cross-sections.
The low number of nodes leads to a very efficient calculation for long parts of the domain.
The outer boundaries parallel to the $x$-direction, illustrated by thin lines in Figure~\ref{fig:Crack}, have traction-free conditions.
Such super elements can also approximate semi-infinite parts of the model that are used to bypass reflections, see, e.g., $\Omega_1$ and $\Omega_6$ in Figure~\ref{fig:Crack}.
For the derivation of the local stiffness matrix of the super element, see \cite{gravenkamp2015simulation,gravenkamp2018efficient}.
The main part of the computation to set up the local stiffness matrix is an eigenvalue problem. A direct formula~\cite{van2007computation} computes the derivative of the eigenvalue problem inside the differentiated algorithm when applying AD.

The super element that forms a star-convex polygon is always required when the prismatic super element can no longer define the geometry or boundary conditions. Here, a finite element line mesh is scaled toward a single point. The single point is called the scaling center. 
In general, the position of the scaling center is also dependent on the parameter set $\para_i$ to be determined by the optimization.
The super element only needs nodes on the boundary yielding a high level of effectiveness. An additional special feature is that the super element also provides a simple crack model because a double node can introduce it. The crack is illustrated by $\Omega_4$ in Figure~\ref{fig:Crack}, where a red node marks the double node. The red line is a traction-free boundary due to the double node. 
Despite the crack tip singularity, one still observes exponential convergence under $p$-refinement~\cite{bulling2019comparison,tschoke2018numerical,bulling2019high} using such super elements.
For the derivation of the local stiffness matrix of the super element, we followed~\cite{chen2014high}. 
The derivation of the local stiffness matrix is quite involved. A continued-fraction-based approach builds on one algebraic Riccati equation and several algebraic Lyapunov equations~\cite{song1997scaled,bazyar2008continued}. For this work, both algebraic equations are solved by an eigenvalue decomposition inside the SBFEM algorithm. The solution for the Riccati equation can be found in~\cite{song1997scaled}, while the Lyapunov equation is solved in the appendix~\ref{sec:ID}.

It is worth mentioning that previous investigations have shown that polygonal super elements are quite robust for small and large finite element sizes inside the same mesh~\cite{xing2018scaled}. This robustness is important if geometric optimization changes the element sizes on the fly.
	\section{The Adapted Optimization Approach}\label{sec:Optimization}

\subsection*{The Inverse Problem}
In general, the procedure described here to identify the defect properties is an indirect method. The aim is to solve the inverse problem
\begin{equation}
\min_{\vpara \in \paraspace} \|\Fo(\vpara) - \vy_{meas}\|^2,
\label{eq:inverse_problem}
\end{equation}
where $\vy_{meas}$ is given by physical measurements performed to detect the defect and
$\Fo(\vpara)\colon \vpara \mapsto \vy_{sim}$ denotes the \textit{Forward Operator}. The forward operator is given by a simulation of the same quantity assuming that the defect is characterized by some vector $\vpara$.
Here, the forward operator is a short notation for the EWM together with the SBFEM-model described in the previous section.
The space $\paraspace$ and its elements $\vpara\in \paraspace$ represent one possible parameterization of the considered defect.
For the examples considered in this article the space  $\paraspace$ is scaled either to $[1,2]^2$ or $[1,2]^3$. This scaling to the same range of values for all parameters supports the solution of the inverse problem. Solving the optimization problem \eqref{eq:inverse_problem}, one obtains a parameterization $\vpara^{\min}$ for the defect which best approximates the measurement data.

In practice, inverse problems are usually hard to solve for several reasons.
Firstly, the resulting optimization problem need not (and generally does not) have a unique solution which may not depend continuously on the measurement data $\vy_{meas}$. 
Such inverse problems are called \textit{ill-posed} problems and are pervailant throughout almost all applications of parameter identification problems.
In order to solve ill-posed problems one is usually well-advised to deploy custom-tailored regularization techniques. 

On the other hand, a different issue that may arise is that the error function in \eqref{eq:inverse_problem} may be hard to optimize as it may have adverse optimization properties such as not being sufficiently smooth or incorporating many only locally optimal points. In this case, it may be advantageous to modify the forward operator $\Fo$ and the measurements $\vy_{meas}$ so that the objective function is better suited for optimization purposes. 

\begin{figure}
\centering
\subcaptionbox[Short Subcaption]{Response and envelope at point $P$ \label{fig:WaterfallCrack1}.}[0.59\textwidth]
{\includegraphics[height=0.24\textheight]{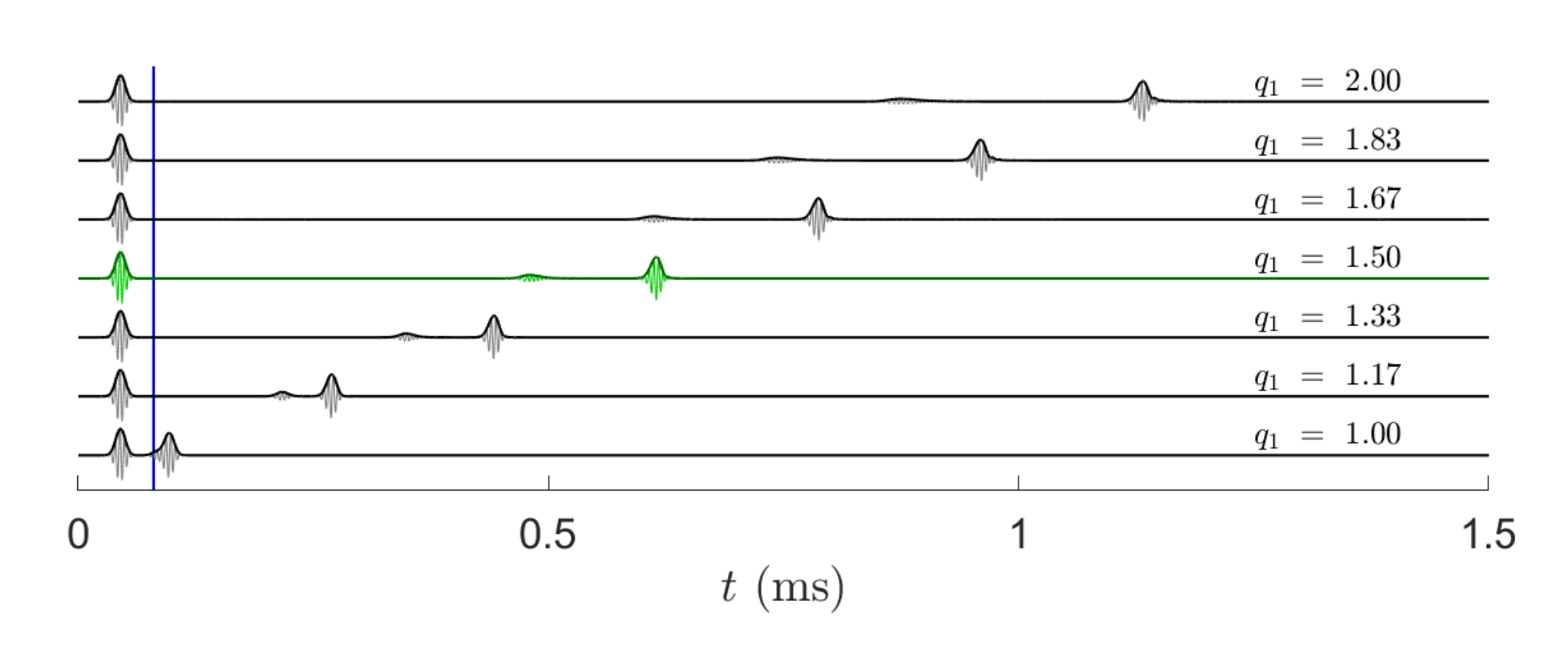}}
\subcaptionbox[Short Subcaption]{Response objective function.\label{fig:WaterfallCrack2}}[0.19\textwidth]
{\includegraphics[height=0.24\textheight]{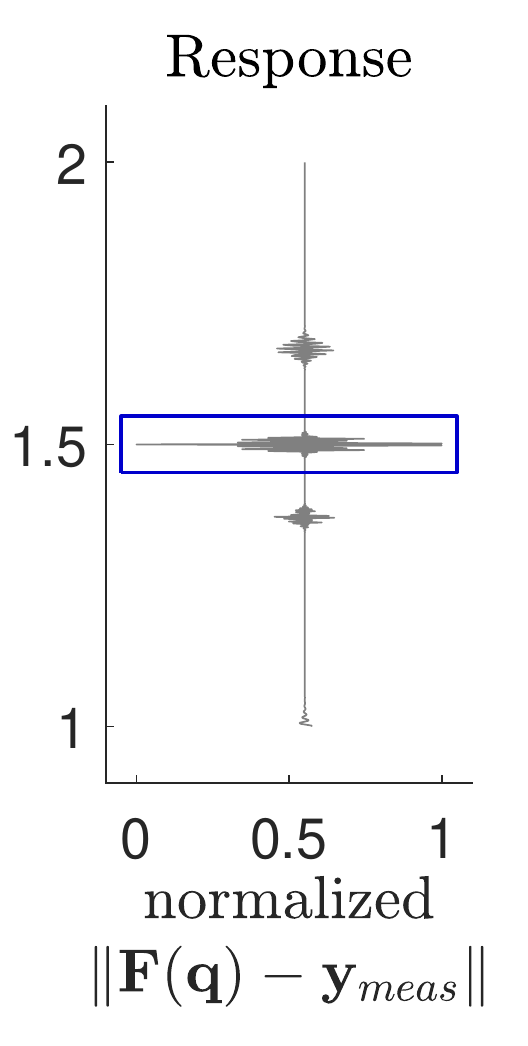}}
\subcaptionbox[Short Subcaption]{Envelope objective function.\label{fig:WaterfallCrack3}}[0.19\textwidth]
{\includegraphics[height=0.24\textheight]{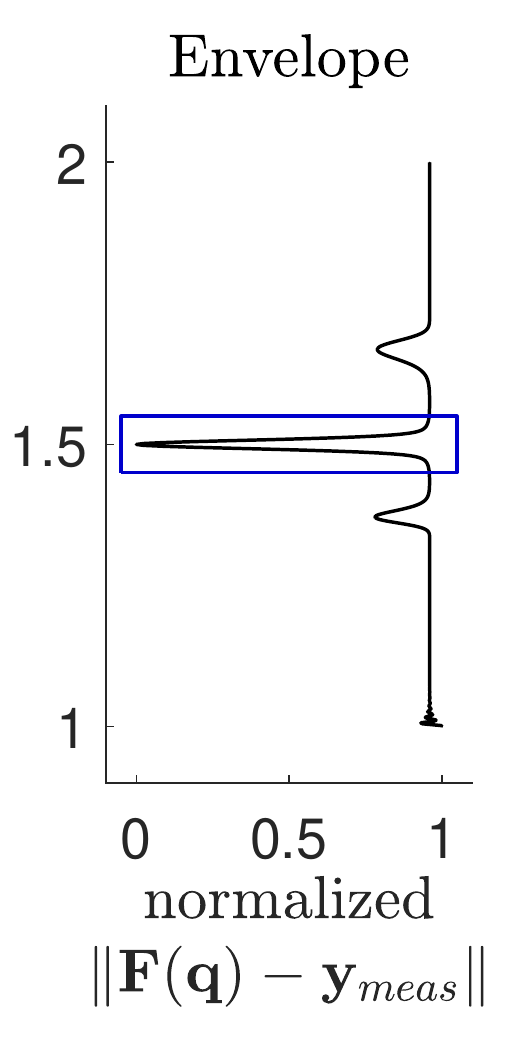}}
\caption{Signals and objective functions for the \textit{waveguide with a crack}.}\label{fig:WaterfallCrack}
\end{figure}

\begin{figure}[htb]
\centering
\subcaptionbox[Short Subcaption]{Response objective\label{Fig:response1}.}[0.45\textwidth ]
{\includegraphics[width=0.43\textwidth]{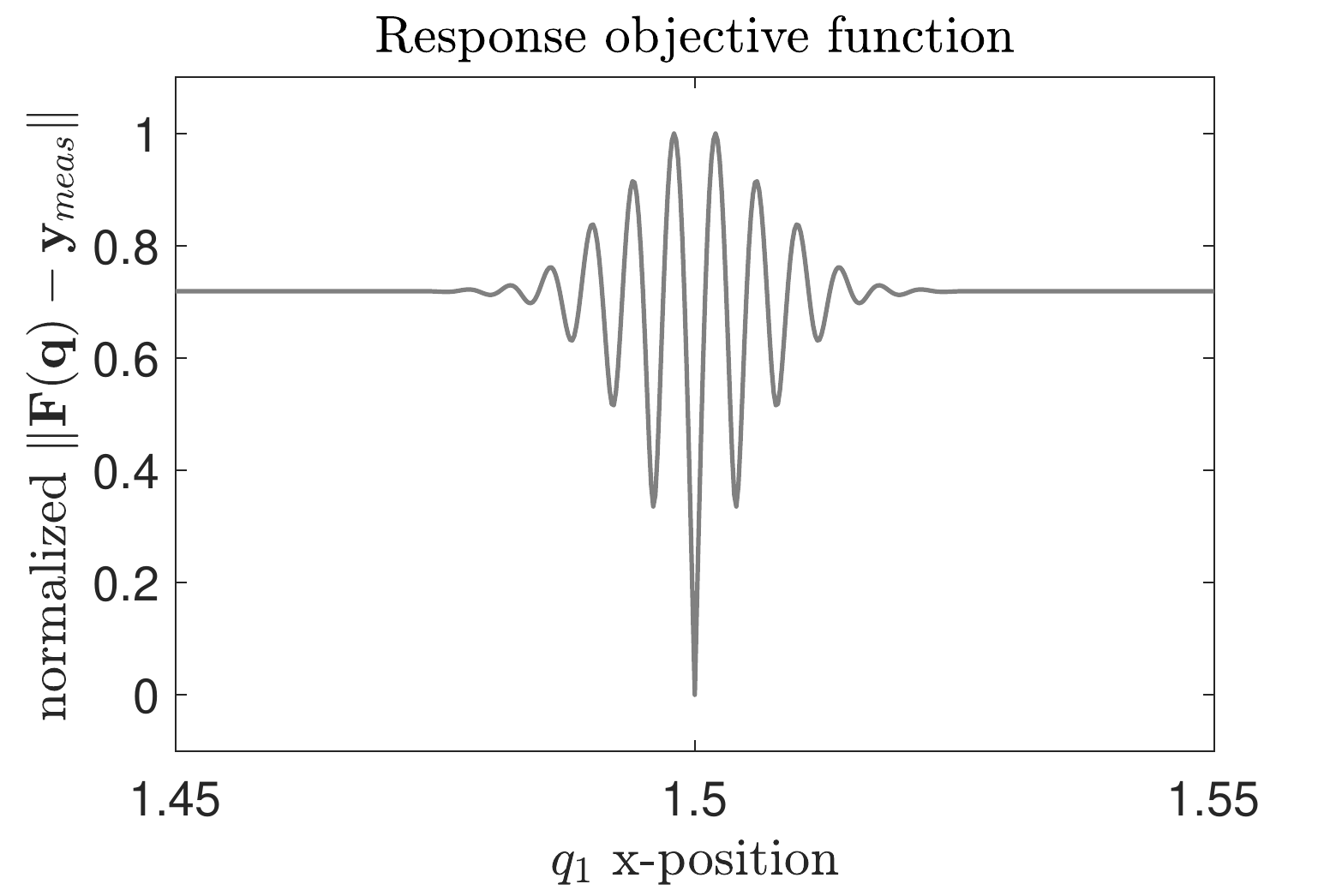}}
\hspace{0.05\textwidth} % seperation
\subcaptionbox[Short Subcaption]{Envelope objective\label{Fig:envelope1}.}[0.45\textwidth]
{\includegraphics[width=0.43\textwidth]{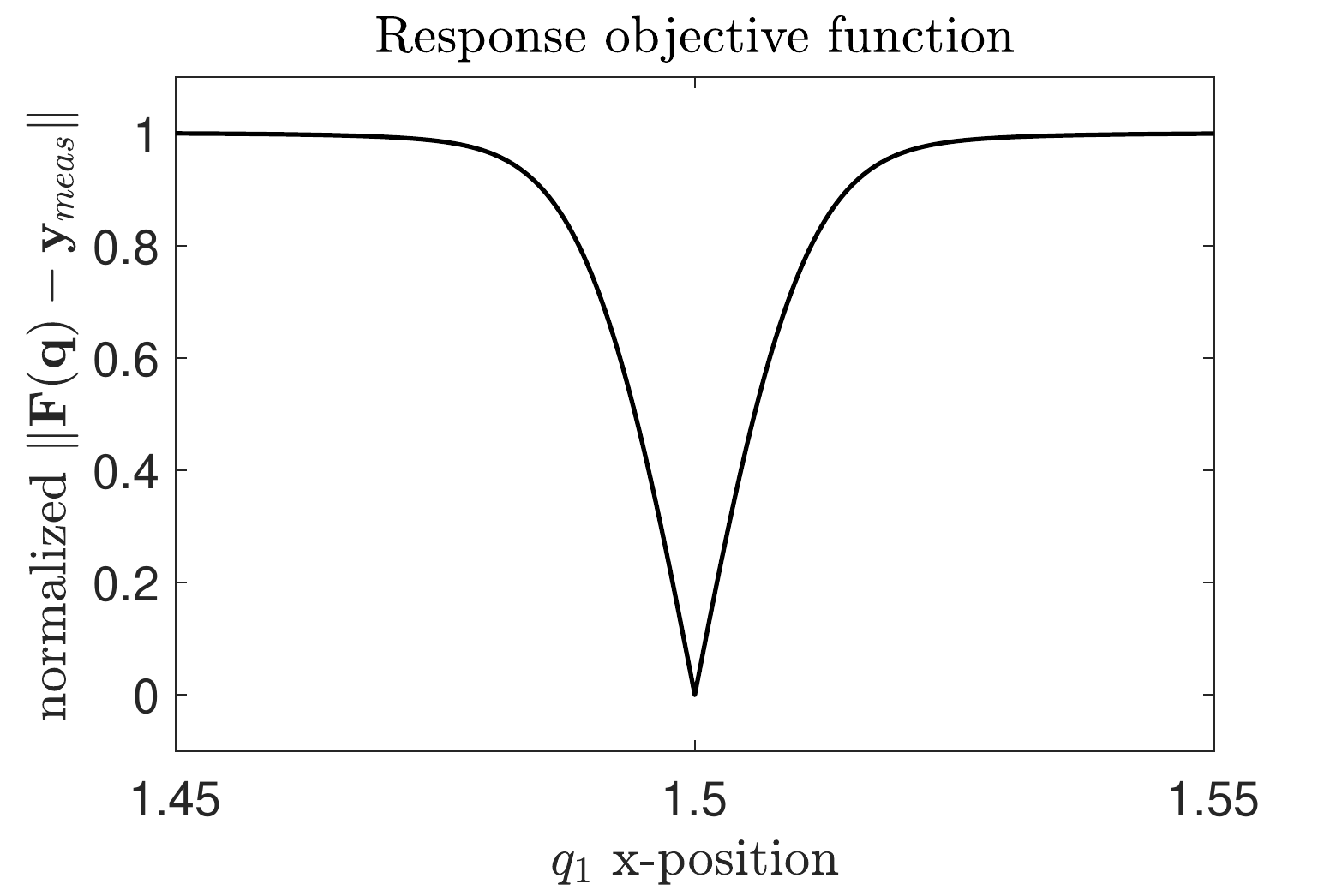}}
\caption[Short Caption]{Comparison of different objective functions as a zoom of the blue box in Figures~\ref{fig:WaterfallCrack2}~and~\ref{fig:WaterfallCrack3}, respectively.}
\label{fig:label}
\end{figure}

A simple approach uses the response data measured directly at a point.
The $y$-component is more interesting as the $x$-component because of its practical relevance. For example, a single laser Doppler vibrometer measures the $y$-component of displacement or velocity, depending on the type of measurement method.
We will assume that the displacement is measured. However, since both quantities have very similar characteristics, all results are transferable to the velocity.
For all numerical examples, the first geometric parameter is the global position of the defect.
Figure~\ref{fig:WaterfallCrack1} shows in gray the $y$-displacement response of the model in Figure~\ref{fig:Crack} at the point $P$ for the global defect $x$-position $\para_1$, where the other parameters, crack depth and angle, are kept constant at $1.5$.
Additionally the envelope of the signal is shown in black or dark green.
Figure~\ref{fig:WaterfallCrack2} depicts the objective function based on the response, where the green signal in Figure~\ref{fig:WaterfallCrack1} is used as measurement data $\vy_{meas} = \Fo(\vpara^*)$.
However, as the waveguide is excited via a sinusoidal pulse, the measured response also incorporates the oscillating behavior.
In fact, it can be seen that the objective function \eqref{eq:inverse_problem} also inherits these features, see Figure~\ref{Fig:response1} for an illustration, which is turned zoom of the blue box in Figure~\ref{fig:WaterfallCrack2}. 
We note that the resulting objective function has many only locally optimal points near the global optima. Hence, many optimization methods would have difficulties solving this optimization problem to global minimality.

The authors also experimented with different variations of the data such as using the spectrum of the response. However, the problem caused by many local minima remained. Many test revealed that the objective function with the best properties from an optimization point of view was the calculation of 
the envelope of the response via a Hilbert-transform, see Figure~\ref{fig:WaterfallCrack3}~and~\ref{Fig:envelope1}. 
Analog to the response, Figure~\ref{fig:WaterfallCrack3} shows the objective function for all crack positions between the green envelope in Figure~\ref{fig:WaterfallCrack2} with a forward operator which included the calculation of the envelope.

In addition to the global minimum, Figure~\ref{fig:WaterfallCrack3} shows two local minima. These local minima are associated with the mode conversion. The excitation, to the left of the blue line, is a 
pulse that corresponds to the S0-mode. However, the S0-mode and the mode converted A0-mode are reflected. If the S0-wavepackage of the forward simulation is at the position of A0-wavepackage of the measurement data, this leads to a local minimum. Compare the green line with the line above in Figure~\ref{fig:WaterfallCrack1}. Our approach to circumvent the local minima is to get a sufficient initial guess, which lies near the global minimum. The details of the initial guess can be found below.

For the rest of the paper, the measurement data is the envelope of the y-displacement at a single point, while the forward operator includes the calculation of the envelope.

%%%%%%%%%%%%%%%%%%%%%%%%%%%%%%%%
%%%%%%%%%%%%%%%%%%%%%%%%%%%%%%%%

%%%%%%%%%%%%%%%%%%%%%%%%%%%%%%%%

\subsection*{Algorithmic Differentiation}
In order to use gradient-based optimization methods for the reconstruction of defects, derivatives of solutions of the partial differential equation are required.
These could be approximated via finite differences (e.g.~\cite{jung2014modeling}) or computed by setting up and solving the adjoint equation of the partial differential equation considered (e.g.~\cite{seidl2016iterative,rao2016guided}).
We tested the use of finite differences and found the approximation to be insufficient for the optimization procedure. 
Furthermore, in the development phase of this project we experimented with different boundary conditions. Hence, an adaptation of the adjoint equation as a time consuming and error-prone process would have been necessary in order to compute derivatives via the adjoint approach. 

Therefore, we apply algorithmic differentiation, also called automatic differentiation (AD), which provides derivative information of arbitrary order for a given code segment. The derivatives are provably accurate up to working accuracy due to the fact that the computer program evaluating the function is broken down into a sequence of elementary evaluations upon which the chain rule of calculus is systematically applied. There are numerous AD tools available for a wide range of programming languages, see the web site www.autodiff.org for an overview of tools and references. 
Since our forward operator described in detail in the previous section is implemented in MATLAB, we  apply the AD-tool ADiMat \cite{adimat} to compute derivatives of solutions of the partial differential equations. 
Some mild modifications to the code were necessary such as computing derivatives of the Lyapunov equation which originally was not available in ADiMat.

\subsection*{The Resulting Optimization Procedure}

Even though we made some effort to reformulate the objective function \eqref{eq:inverse_problem} so that it is easier to optimize, we cannot eliminate all issues entirely. 
Hence, in this subsection, we present the optimization algorithms that we deployed.

A first simple approach for the solution of the optimization problem \eqref{eq:inverse_problem} is to use standard optimization software naively.
This can involve derivative-based approaches such as \texttt{IPOPT} \cite{ipopt}, \texttt{WORHP} \cite{worhp}, or the optimization routine \texttt{fmincon} provided by Matlab. Alternatively, one may
use derivative-free approaches such as coordinate descent,  genetic optimization or other heuristic-based optimization methods like particle-swarm methods.

Derivative-based optimization algorithms have the great benefit of fast and provable convergence but generally only converge locally. With derivative-free methods often there is the hope associated that these methods converge globally.  However, global convergence cannot be proven. Moreover, the whole optimization procedure may require thousands of function evaluations. 
The required number of evaluations is of particular significance as the simulation in our case involves solving a partial differential equation; thus one objective function evaluation is already quite expensive -- solving the partial differential equations many times is simply impractical.

In light of these facts, the approach presented here is a two-step approach:
We first generate a refined initial guess for all variables.
After a good initial point is established, we deploy a gradient-based optimization with fast convergence properties.

\paragraph{Initial guess for the global position} As mentioned above, the first geometric parameter $\para_1$ is the global position of the defect for all numerical examples. To get an initial guess for the defect position $\para_1$, we consider the time of flight to get a physical position $\Ldefect$ inside the waveguide. This physical position is then mapped linearly into the parameter space. The physical position of the defect $\Ldefect$ can be estimated by
\begin{equation}
	\Ldefect \approx \Tdiff \frac{\cgP \cdot \cgN}{\cgP + \cgN},
	\label{eq:initguess_q1}
\end{equation}
where $\Tdiff$ is the time of flight, $\cgP$ is the group velocity at the center frequency of the excited mode, and $\cgN$ is the group velocity of the reflected mode with the highest amplitude.
For all our examples prior investigations have shown that the highest reflections for the $y$-component are always associated with the A0-mode (see Fig.~\ref{fig:WaterfallCrack1}), so the reflection velocity $\cgN$ is the A0-group velocity $\cg^\text{A0}$ at the center frequency - see Figure~\ref{fig:excitationB}.

The computation of the time of flight $\Tdiff$ is based on cross-correlating the excitation pulse, which is separated by a defined threshold $\Tex$ from the reflections, with the rest of the signal.
If the envelope of the cross-correlation is maximal, the argument is an approximation of the time of flight, i.e.,
\begin{align}
    r & = \env\left( u_y(t)[t<\Tex] \star u_y(t)[\Tex \leq t] \right),\\
    \Tdiff & = {\arg \max}_t r,
\end{align}
where $u_y$ is the $y$-displacement, $\env$ is the envelope function, $\star$ is the cross-correlation operator and $[\cdot]$ is the Iverson bracket\footnote{The Iverson bracket is one if the statement inside is true and zero otherwise.}.
Figure~\ref{fig:WaterfallCrack1} shows the threshold as a blue line, and the excitation pulse is on the left side while the rest of the signal is on the right side.
In the presented examples, the function $r$ has a distinct maximum, which is associated with the A0-mode. However, if there is any ambiguity in the mode or the maximum value, it is simple to run the inverse method with multiple start values for the first parameter.

\paragraph*{Refined initial guess for all parameters}

In order to obtain a refined initial guess for the gradient-based optimization, we can use the approximation Equation~\eqref{eq:initguess_q1}.
The approximation for the first coordinate given there is fairly accurate and is located within the valley of the objective function, see Figure \ref{Fig:envelope1}.
However, with a simple line-search in the vicinity of the valley, the accuracy of the initial guess can be greatly improved. 
This simple line-search requires very few function evaluations, and the computational cost is insignificant.
The accuracy of the first coordinate is important as inaccurate values lead to locally optimal values (especially in the presence of noise) for the other coordinates with which optimization methods have difficulty escaping. 
The other coordinates are determined by selecting the vector with the smallest objective function from $\Nsample$ random vectors, where our default value of $\Nsample$ is 10.
As the line-search is computationally very cheap, we also refine the other coordinates in the same fashion, beginning with the least critical or problematic coordinates.
Due to many only locally optimal points, see Figure~\ref{Fig:Corrosion_objective_function}, the amount of randomly selected points is increased to 100 for corrosion-type defects.
The whole procedure is sketched in Algorithm~\ref{Algo1}.

\begin{algorithm}
\caption{Refined initial guess}
\begin{algorithmic}[1]
 \State Evaluate Equation~\eqref{eq:initguess_q1} for an initial guess of the first coordinate.
 \State Sample $\Nsample$ random points $\vpara_0^1,\dots,\vpara_0^\Nsample \in \paraspace$ for the remaining second and possibly third coordinate.
 \State Evaluate corresponding objective values $\|\Fo(\vpara_0^m) - \vy_{meas}\|^2$, $m = 1,\ldots,\Nsample$, and select $\vpara_0^m$ with the smallest objective value. 
 \State Apply coordinate-wise line-search to $\vpara_0^m$ to decrease the objective function value. The resulting parameter is $\vpara_0$.
\end{algorithmic}
\label{Algo1}
\end{algorithm}

\subsubsection*{Optimization procedure}
We experimented with different optimization methods in order to solve the optimization problem Equation~\eqref{eq:inverse_problem}. 
One difficulty we encountered was that the objective function has many locally optimal points that are difficult to escape from. 
However, optimization approaches aiming at a more global scale, such as genetic algorithms, or particle swarm methods are infeasible due to the computation effort caused by one evaluation of the forward operator.

We were able overcome these issues by applying an Iteratively Regularized Gauss-Newton Method (IRGNM), see e.g. \cite{kaltenbacherbuch}.
This method modifies the well known Gauss-Newton approach for nonlinear least-squares fitting problems by adding terms that make sure that matrices occurring in the Newton-step are invertible and adds a Thikonov-type term so that the overall solution stays near the initial guess.
The overall procedure is given in Algorithm~\ref{Algo2}, where $\Fo'$ is the Jacobian matrix of the forward operator, $\mI$ is the Identity matrix and $(\cdot)^\ct$ denotes the adjoint operator.

\begin{algorithm}
\caption{Iteratively Regularized Gauss-Newton Method (IRGNM)}
\begin{algorithmic}[1]
 \State Initialize $\vpara_0$ and  $(\alpha_n)_{n\in {\NN}}\rightarrow 0$
  \For{$n = 0,\dots, {maxiter}$}
 \State $\vpara_{n+1} = \vpara_n + (\Fo'(\vpara_n)^\ct\Fo'(\vpara_n)+ \alpha_n \mI)^{-1} (\Fo'(\vpara_n)^\ct (\vy_{meas}-\Fo(\vpara_n)) + \alpha_n (\vpara_0-\vpara_n) )$
 \If{$\|\vpara_{n+1}-\vpara_n\|<\epsilon$,}
    \State $\vpara^{\min} = \vpara_{n+1}$
    \State STOP
 \EndIf
 \EndFor
\end{algorithmic}
\label{Algo2}
\end{algorithm}

	\section{Numerical Experiments}\label{sec:NumericalExperiments}

In this section, we present several numerical examples to show the versatility of the approach. Each model considers another type of defect. The first example is about classic crack identification. The second example shows a similar set-up but tries to find a horizontal crack or rather a delamination defect. The third type of defect is a model of a corrosion defect.

All waveguides in the following examples are made out of steel. The isotropic material parameters are given in Table~\ref{tab:Steel}, and for all waveguides a plane strain is assumed. All meshes use shape functions of degree $\pp = 4$, i.e., elements with five nodes. The necessary degree is determined in advance by investigating the critical geometric parameter. For example, the smallest and the largest crack length are the critical geometric parameter for the first model. For all our models, a simulation with a degree of $\pp = 4$ leads to a difference in the signals below $0.1\%$ compared to a simulation with $\pp = 5$. It is worth noting that there are numerous convergence studies for the SBFEM~\cite{tschoke2018numerical,gravenkamp2020mass,bulling2019comparison,chen2014high} and the code is thoroughly validated using commercial FEM tools.

\begin{table}[htb]
    \centering
    \begin{tabular}
    {rS[table-format = 3.0]s|
     rS[table-format = 1.1]s|
     rS[table-format = 1.2]s}
         \multicolumn{9}{c}{Isotropic material} \\
         \toprule
         $E\colon$ & 200&\giga\pascal&
         $\nu\colon$& 0.3 & &
         $\rho\colon$ & 7.85&\g\per\cubic\cm\\
    \end{tabular}
    \caption{Structural steel}
    \label{tab:Steel}
\end{table}

\subsection{Waveguide with a Crack}\label{sec:Crack}

\begin{table}[htb]
    \centering
    \begin{tabular}
    {r|cS[table-format = 4.1]S[table-format = 4.1]s}
    	{description}&{parameter}&{min}&{max}&{Unit}\\ \hline
		global $x$-position of $\Omega_5$&$\para_1$&200&2200&\milli\meter \\
	    $x$-position of the crack tip with 0 as the midpoint of $\Omega_5$&$\para_2$&-0.5&+0.5&\milli\meter \\
		y-position of the crack tip with 0 as the midpoint of $\Omega_5$&$\para_3$&-2.25&0&\milli\meter \\
    \end{tabular}
    \caption{Parameters for the \textit{waveguide with a crack}.} \label{tab:CrackPara}
\end{table}

\begin{figure}
\centering
\subcaptionbox[Short Subcaption]{Objective function.\label{fig:cracked_plate_obj}}[0.45\textwidth ]
{\includegraphics[width=0.43\textwidth]{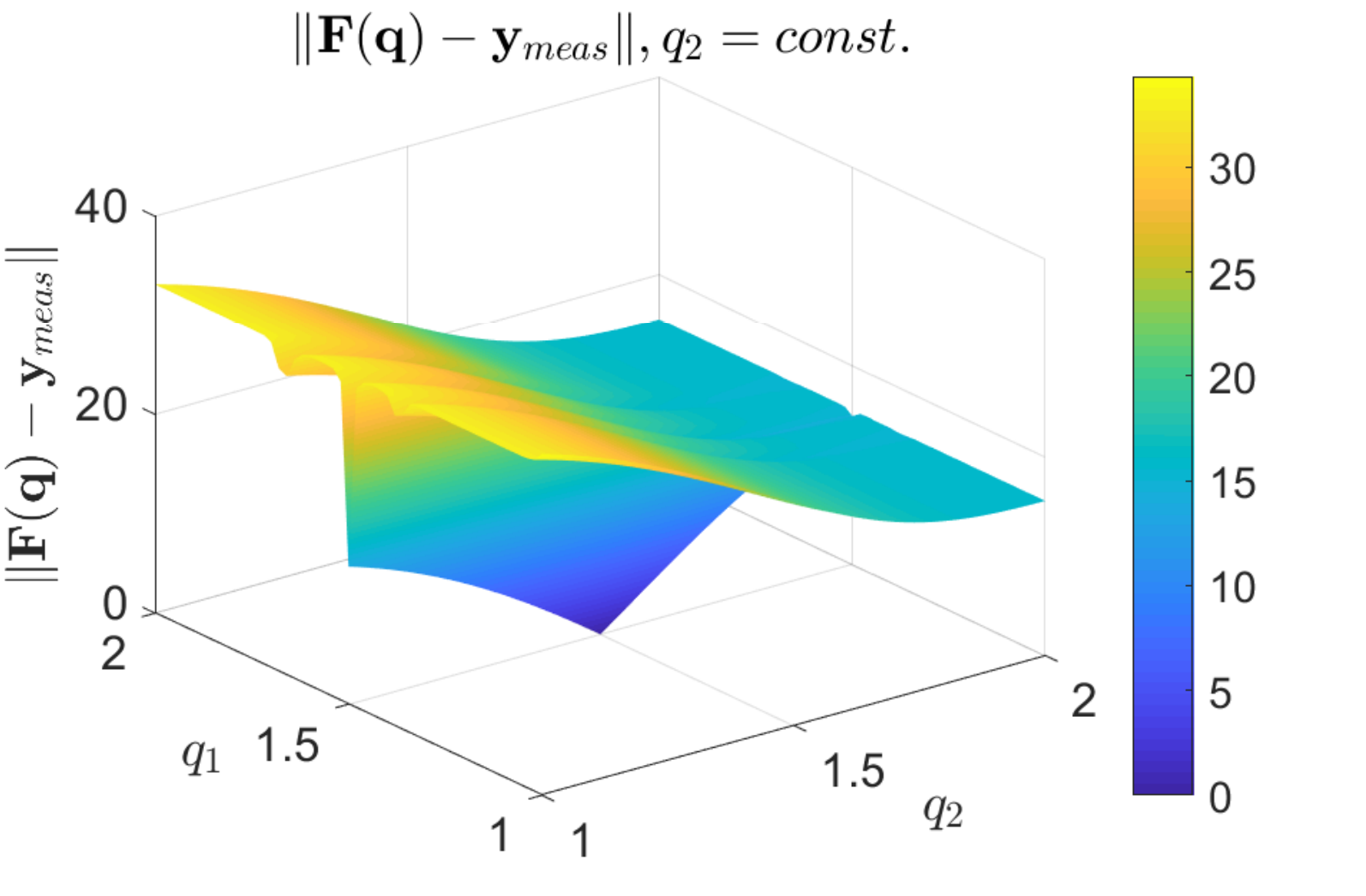} }
\hspace{0.05\textwidth} % seperation
\subcaptionbox[Short Subcaption]{Reconstruction error of parameters\\ $|\vpara_i^{\min}-\vpara_i^*|$ with noisy data. \label{fig:crack_plate_noise}}[0.45\textwidth]
{\includegraphics[width=0.43\textwidth]{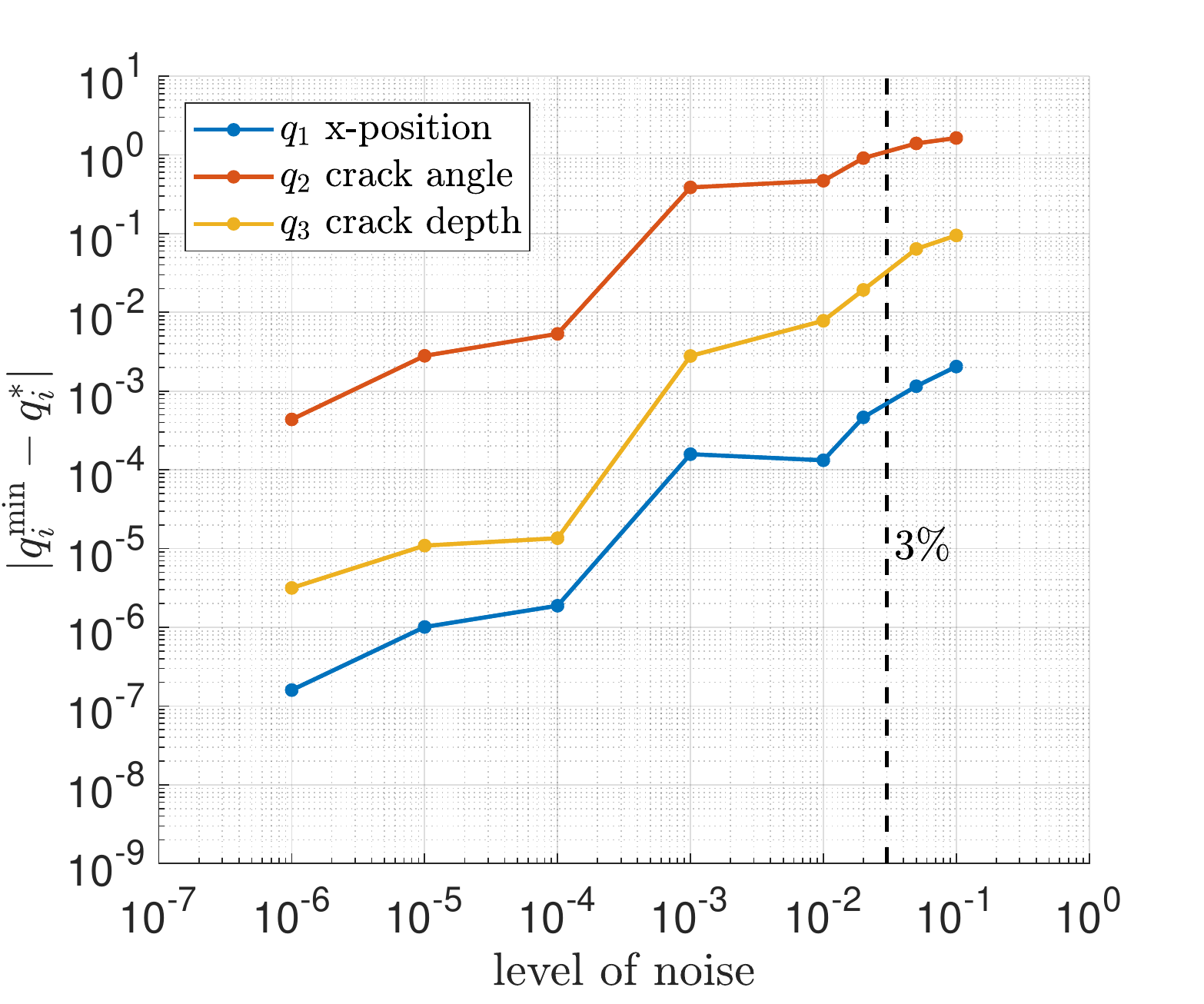} }
\caption[Short Caption]{The objective function and reconstruction with noisy data for the \textit{waveguide with a crack}, where the target parameter $\vpara^*$ and the associated artificial measurement signal $\vy_{meas}$ correspond to the midpoint of the parameter space.}
\end{figure}

\begin{figure}
\centering
\subcaptionbox[Short Subcaption]{\textit{Waveguide with a crack}. \label{Fig:iterates_crack}}[0.45\textwidth]
{\includegraphics[width=0.43\textwidth]{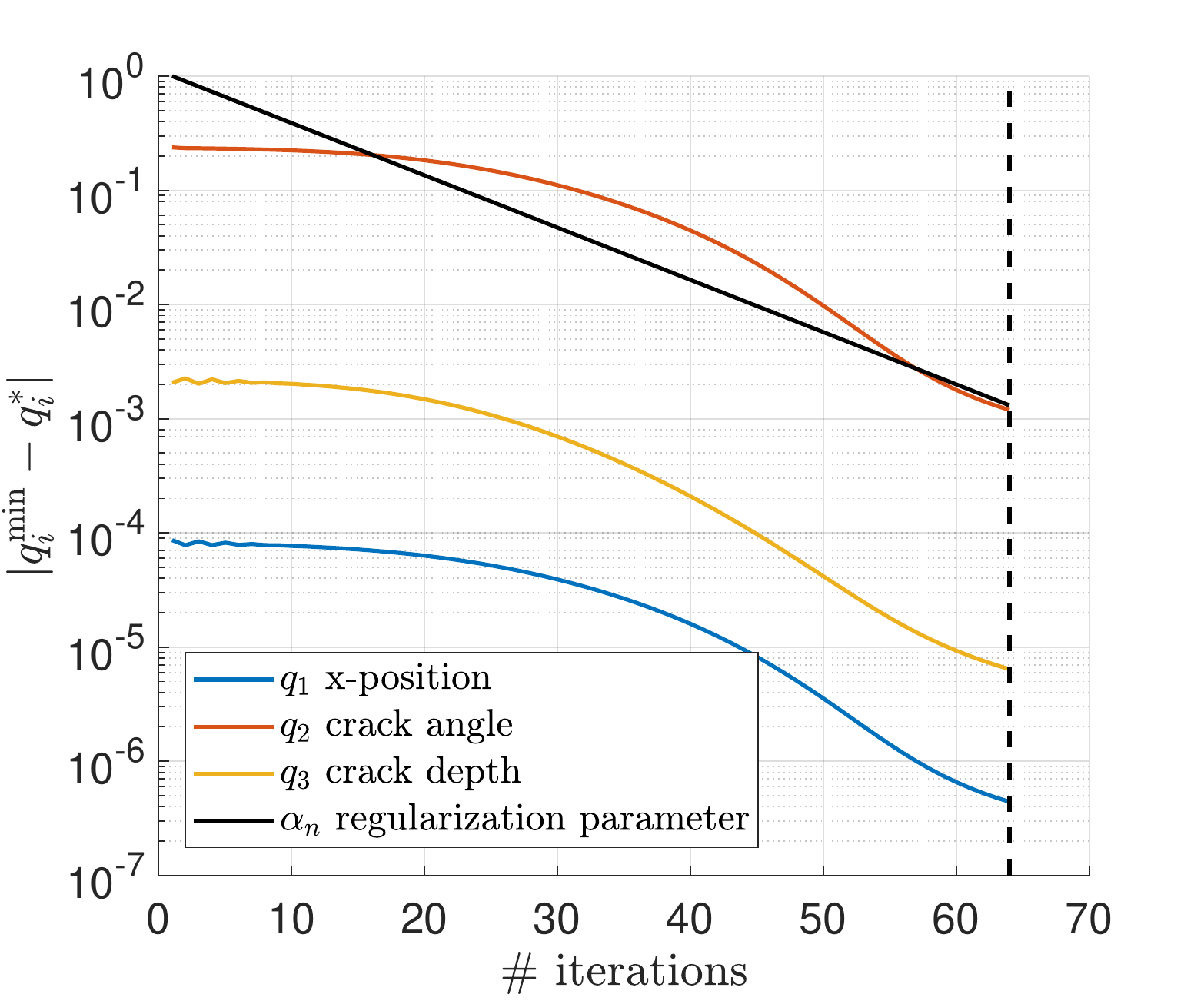}}
\hspace{0.05\textwidth} % seperation
\subcaptionbox[Short Subcaption]{\textit{Waveguide with a corrosion defect}. \label{Fig:iterates_corro}}[0.45\textwidth]
{\includegraphics[width=0.43\textwidth]{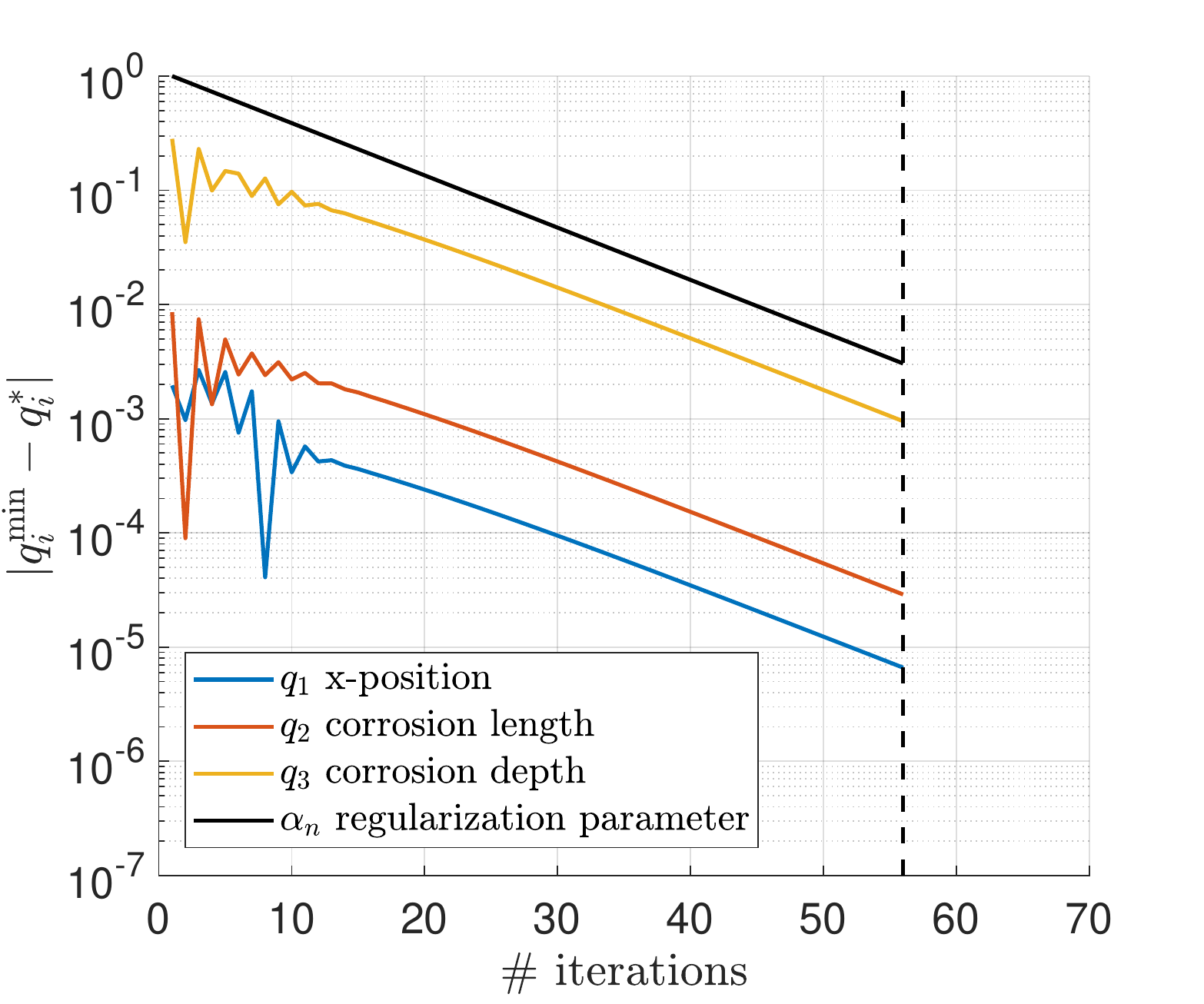}}
\caption[Short Caption]{Error in reconstruction procedure, where the target parameter $\vpara^*$ is the midpoint of the parameter space and  a small noise of $10^{-5}$ is added.}
\end{figure}

\begin{table}[htb]
\centering
\begin{tabular}{r|llllllllll}
$\para_1^*$& 1.33& 1.73& 1.72& 1.58& 1.58& 1.43& 1.46& 1.54& 1.71& 1.31 \\
$\para_2^*$& 1.25& 1.54& 1.57& 1.37& 1.50& 1.52& 1.34& 1.42& 1.40& 1.37 \\
$\para_3^*$& 1.47& 1.67& 1.74& 1.69& 1.68& 1.38& 1.41& 1.40& 1.47& 1.46 \\ \hline
$\|\vpara^{\min}-\vpara^*\|$& 2E-5& 3E-3& 2E-3& 6E-4& 5E-3& 7E-4& 1E-4& 1E-4& 4E-4& 5E-6
\end{tabular}
\caption{Error for reconstruction $\vpara^{\min}$ for differently shaped defects $\vpara^*\in [1,2]^3$ for the \textit{waveguide with a crack} with a noise level of $10^{-5}$.}
\label{tab:errors_crack}
\end{table}

This example considers a steel plate with a single straight crack. Figure~\ref{fig:Crack} shows an overview of the model. Note that the $x$-axis is interrupted to fit the model into the figure. The first section of this figure is for the excitation. Arrows in Figure~\ref{fig:Crack} show the area and direction where spatially constant traction is applied. At these parts, the excitation takes place by a normal traction force
\begin{equation}
    \ff_{S0}(\vx) = - \vn [\vx \in \Gamma_\tau] \label{eq:tractionS0}.
\end{equation}
The traction is a model for the excitation of a double transducer set-up~\cite{su2004selective}. A double transducer set-up can produce both fundamental modes by applying a force in the same or opposite direction on the plate.
For this traction, the S0-mode is excited, so the group velocity $\cgP$ in equation~\eqref{eq:initguess_q1} is S0-group velocity $\cg^\text{S0}$. Here, the S0-excitation is appropriate because the S0-mode leads to larger reflections for these kinds of cracks.

The second section of Figure~\ref{fig:Crack} shows the evaluation point $P$. The envelope of the $y$-displacement is evaluated as a time series. As described above, the $y$-component is of practical relevance as it can be measured using a single laser Doppler vibrometer.
Figure~\ref{fig:WaterfallCrack} shows the $y$-displacement at point $P$ in gray and the envelope as a black curve over the signal of different crack positions. 

The third section contains the defect. A double node is inserted into a continued-fraction-based super element. The emerging inner traction-free boundary is marked in red in Figure~\ref{fig:Crack}. Note that the continued-fraction approach handles the singular stress at the end of this crack model. In total, the model has $180$ degrees of freedom. This low number is possible because the semi-analytical approach approximates the large undamaged part of the waveguide. At this center frequency, the wavelength for the S0 mode is around \SI{25}{\milli\meter}. With other simulation methods, this leads to a much higher computational effort.

There are three parameters for the optimization. These parameters change the defect and its location. The first parameter $\para_1$ defines the global position of the crack. This change in position is achieved by shifting the super element $\Omega_5$ along the $x$-axis and changing the length of the super element $\Omega_4$ accordingly. The second and third parameters control the crack length and angle by moving the crack tip inside the super element $\Omega_5$ relative to the boundary. The maximal and minimal values for the parameters are listed in Table~\ref{tab:CrackPara}. The crack length is between \SI{0.25}{\milli\meter} and ca. \SI{2.6}{\milli\meter}.
Here the parameters are restricted to allow a numerically stable solution. Similar parameter restrictions can be found in other publications like \cite{agathos2018multiple}.

The objective function dependent on $q_1$ and $q_2$ with regard to the optimization problem \eqref{eq:inverse_problem} is displayed in Fig.~\ref{fig:cracked_plate_obj}. Here, $q_3$ is kept constant at $1.5$, i.e., only orthogonal cracks are considered.
The artificial target measurement corresponds to the midpoint of the parameter space. The surface is relatively smooth and has only a few local minima. The global minimum at $(1.5,1.5,1.5)^\tr$ has clearly the lowest value.
Using the result of Algorithm \ref{Algo1} as refined initial guess, the defect properties can be obtained by applying the IRGN method Algorithm \ref{Algo2} without any need for using a more in-depth initial guess search such as applying genetic algorithms.

Another important aspect of an inverse method, especially in view of the measurement data, is its performance under noise. To analyze the performance, we create artificial measurement data by a single forward simulation with a parameter vector $\vpara^*$ and add random noise to this target signal. The noise is defined in terms of the excitation pulse. The same threshold as for the time of flight is used to identify the excitation pulse. 
In Figure~\ref{fig:WaterfallCrack}, the blue line shows this threshold. The excitation pulse is to the left of the blue line. A noise level of $10^{- 2}$ and $10^{- 3}$ means that a normal-distributed noise with vanishing mean value a standard deviation of $1\%$ and $0.1\%$ of the maximum amplitude of the excitation pulse has been added to the target signal, respectively. The envelope of the target signal is calculated afterward. Figure~\ref{fig:crack_plate_noise} shows the reconstruction error $|\para_i^{\min} - \para_{i}^*|$ per parameter for the midpoint of the parameter space.
Here, the first parameter, the global position of the crack, shows the smallest error. A lower error in the first parameter is also expected because the range in the physical space is much larger compared to the other parameters – see Table~\ref{tab:CrackPara}. Translating the error in the physical space, the error stays under \SI{1}{\milli\meter} for a noise level under $3\%$, which is a satisfactory level of precision for most applications. The error in the second parameter, on the other hand, shows lower robustness to noise. This parameter, which is linked to the crack angle, is likely to produce incorrect results. However, the more critical third parameter, which is more related to the crack length, has a significantly lower reconstruction error. In general, this investigation is only valid for the current target parameter $\vpara^*$; as the crack length gets smaller with the third parameter, it is also expected that error increases due to the negative influence of the noise. In conclusion, there are also physical boundaries to an appropriate parameter space depending on the noise inside the experimental data. 

Figure~\ref{Fig:iterates_crack} shows the convergence of the reconstruction error with a small noise level of $10^{-5}$. This noise level is introduced to prevent an "Inverse Crime" as both the target data and the reconstruction algorithm use the same simulation method. We observe an asymptotic rate that coincides with the rate of the regulation parameter in the Algorithm~\ref{Algo2}.

To further validate the inverse method, ten different target parameters are randomly drawn. Noise with a standard deviation of $10^{-5}$ is added to the corresponding signals, and the parameters are reconstructed. Table~\ref{tab:errors_crack} lists the reconstruction error for these points. For all reconstructions, the error is below $1\%$.
%%%%%%%%%%%%%%%%%%%%%%%%%%%%%%%%%%%%%%%%%%%%%%%%%%
\subsection{Waveguide with a Delamination}

\begin{figure}[ht]
    \begin{center}
    \includegraphics[width=1\textwidth]{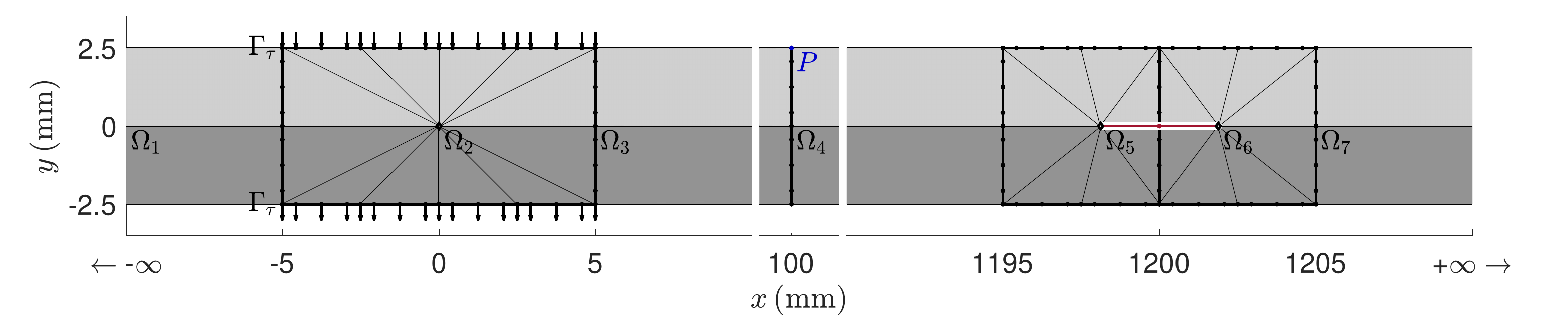}
	\caption{Domain for the \textit{waveguide with a delamination}.}
	\label{fig:Delam}
    \end{center}
\end{figure}

\begin{figure}
    \centering
    \subcaptionbox[Short Subcaption]{Objective function.\label{fig:Delam_obj}}[0.45\textwidth ]
    {\includegraphics[width=0.43\textwidth]{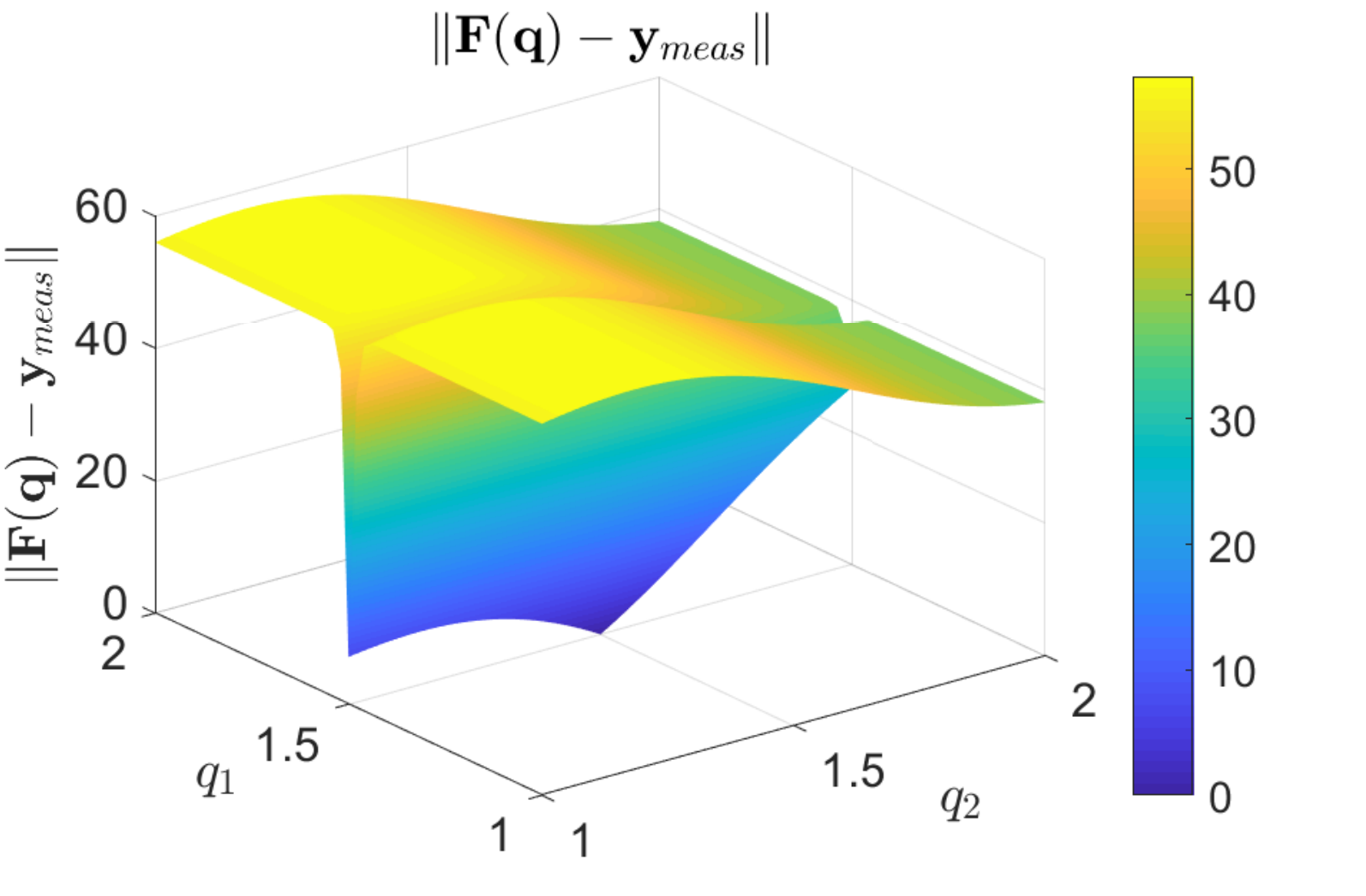} }
    \hspace{0.05\textwidth} % seperation
    \subcaptionbox[Short Subcaption]{Reconstruction error of parameters\\ $|\vpara_i^{\min}-\vpara_i^*|$ with noisy data. \label{fig:delam_noisy_reconstruction}}[0.45\textwidth]
    {\includegraphics[width=0.43\textwidth]{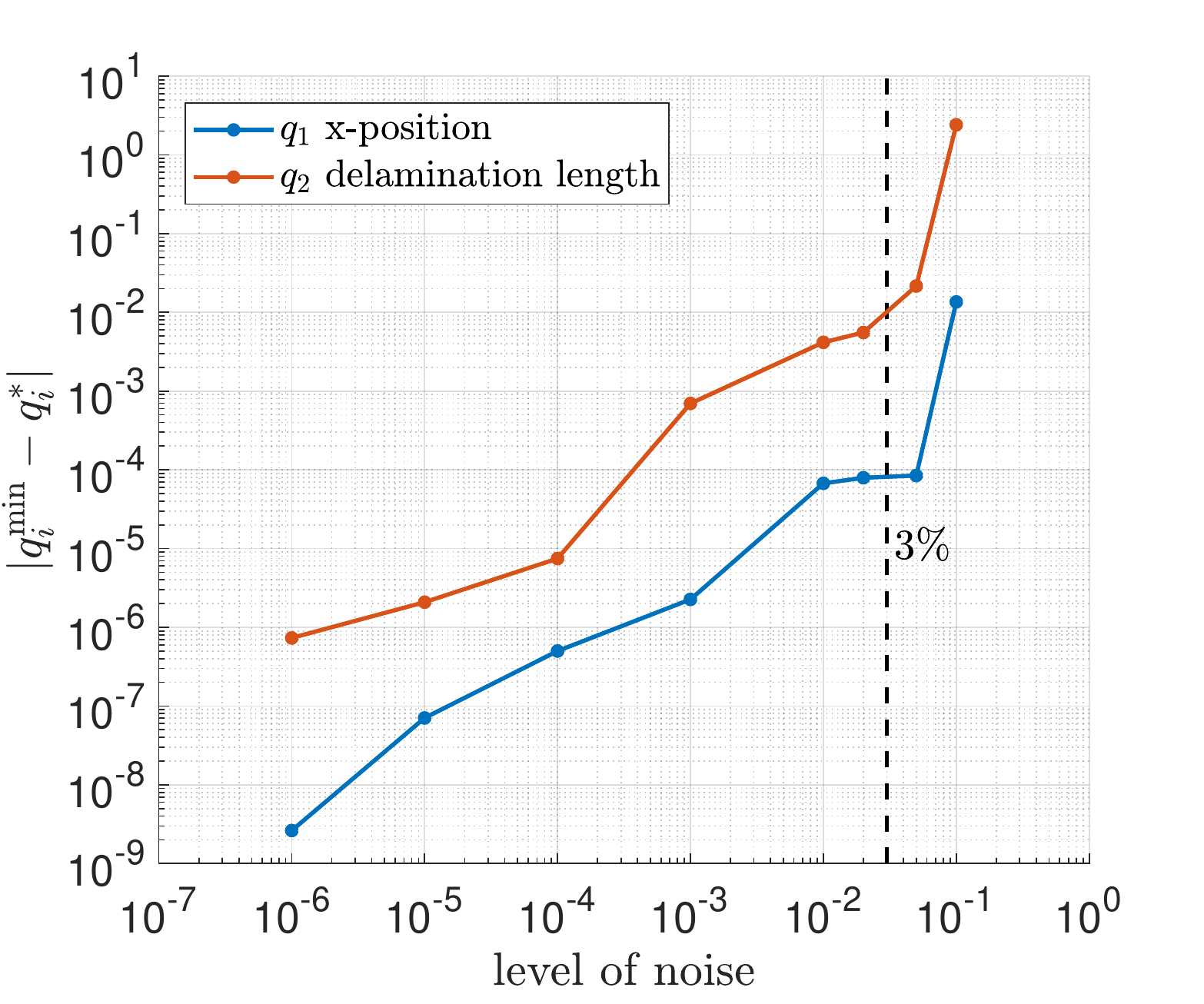} }
    \caption[Short Caption]{
    The objective function and reconstruction with noisy data for the \textit{waveguide with a delamination}, where the target parameter $\vpara^*$ and the associated artificial measurement signal $\vy_{meas}$ correspond to the midpoint of the parameter space.}
\end{figure}

\begin{table}[htb]
    \centering
    \begin{tabular}
    {r|cS[table-format = 4.1]S[table-format = 4.1]s}
    	{description}&{parameter}&{min}&{max}&{Unit}\\ \hline
		global $x$-position of $\Omega_5,\Omega_6$&$\para_1$&200&2200&\milli\meter \\
	    delamination length inside $\Omega_5,\Omega_6$&$\para_2$&2.5&7.5&\milli\meter \\
    \end{tabular}
    \caption{Parameters for the \textit{waveguide with a delamination}.}
    \label{tab:DelamPara}
\end{table}

\begin{table}
    \centering
    \begin{tabular}{r|llllllllll}
    $\para_1^*$& 1.45& 1.28& 1.29& 1.35& 1.66& 1.54& 1.33& 1.47& 1.26& 1.69 \\
    $\para_2^*$& 1.69& 1.38& 1.27& 1.59& 1.67& 1.39& 1.44& 1.43& 1.29& 1.32 \\ \hline
    $\|\vpara^{\min}-\vpara^*\|$& 1E-6& 1E-5& 7E-6& 2E-7& 2E-6& 8E-6& 4E-6& 3E-6& 1E-5& 9E-6
    \end{tabular}
    \caption{Error for reconstruction $\vpara^{\min}$ for differently shaped defects $\vpara^*\in [1,2]^2$ for the \textit{waveguide with a delamination} with a noise level of $10^{-5}$.}
    \label{tab:errors_delam}
\end{table}

This example considers a  plate with a single straight horizontal delamination model.
Figure~\ref{fig:Delam} shows an overview of the model.
Note that once again, the $x$-axis is interrupted at some places to fit the model into the figure.
The two colors inside the figure indicate that these defects are often present between two different layers. However, the material of these layers in our example is the same steel.

Prior investigations have shown that the A0-excitation will lead to larger reflections for delaminations. The traction is chosen accordingly as
\begin{equation}
    \ff_{A0}(\vx) = - \ve_2 [\vx \in \Gamma_\tau],
\end{equation}
where $\ve_2$ is the second unit vector. Hence, $\cgP$ in  equation~\eqref{eq:initguess_q1} is the A0-group velocity $\cg^\text{A0}$ at the center frequency - see Figure~\ref{fig:excitationB}.

The second section in Figure~\ref{fig:Delam}
highlights again the evaluation point $P$, where the envelope of the $y$-component of the displacement is taken as the data.

The third section of Fig.~\ref{fig:Delam} contains the defect. The delamination is modeled by inserting double nodes in two super elements.
The semi-analytical solution of the SBFEM once again captures the points with singular stress.

There are two parameters for the optimization. The first parameter $\para_1$ denotes the global position of the delamination. The second parameter changes the length of the delamination.
The maximal and minimal values for the parameters are listed in Table~\ref{tab:DelamPara}.

Figure~\ref{fig:Delam_obj} depicts the objective function dependent on $q_1$ and $q_2$.
As in the previous example, the target signal is generated with the midpoint of the parameter space for $\vpara^*$.
Figure~\ref{fig:Delam_obj} and Figure~\ref{fig:cracked_plate_obj} show a similar shape, but the two side channels with local minima are not present in this example. These side channels are missing because delaminations in the center of a waveguide do not lead to mode conversion for the A0-mode.

For this example, we could set the regularization parameter $\alpha_n=0$ and achieve convergence within five iteration steps.
In Figure \ref{fig:delam_noisy_reconstruction}, the error for reconstruction with noise is given. Both parameters show high robustness against the noise. Additionally, Table~\ref{tab:errors_delam} lists the reconstruction error for ten different randomly drawn target parameters with a noise level of $10^{-5}$. The error indicates a high level of precision.

In total, the model has $226$ degrees of freedom. A single forward calculation is performed in under $4$ seconds on a modern desktop computer, while the reconstruction is computed in under $5$ minutes, including all preprocessing steps. The time is averaged over 100 reconstructions.
%\TODO{1 Auswertung: 3.71sec (avg von 100 auswertungen), reconstruction 214sec}
%%%%%%%%%%%%%%%%%%%%%%%%%%%%%%%%%%%%%
\subsection{Waveguide with a Corrosion Defect}

\begin{figure}[ht]
    \begin{center}
    \includegraphics[width=1\textwidth]{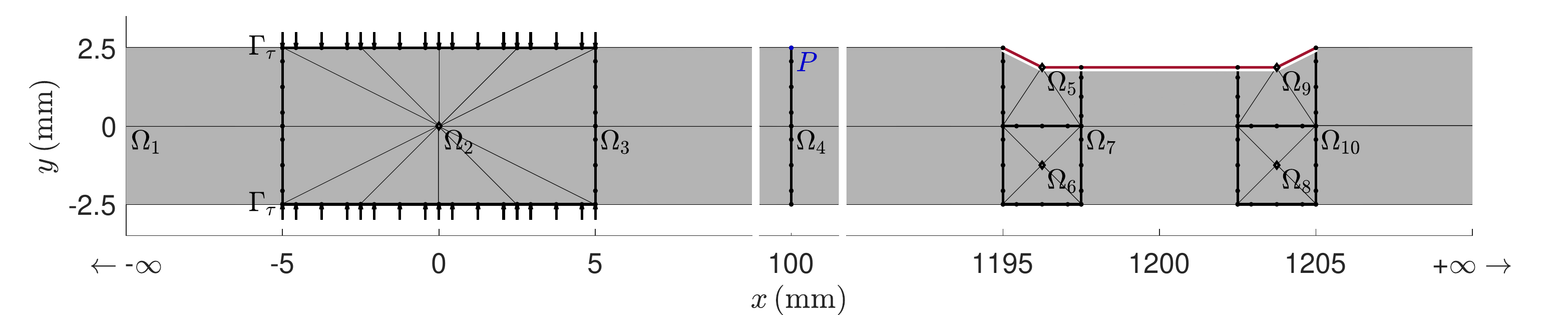}
		\caption{Domain for the \textit{waveguide with a corrosion defect}.}
		\label{fig:Corro}
    \end{center}
\end{figure}

\begin{table}[htb]
    \centering
    \begin{tabular}
    {r|cS[table-format = 3.1]S[table-format = 3.1]s}
    	{description}&{parameter}&{min}&{max}&{Unit}\\ \hline
		global $x$-position of $\Omega_5,\ldots,\Omega_9$&$\para_1$&200&2200&\milli\meter \\
	    corrosion length on the surface&$\para_2$&7.5&57.5&\milli\meter \\
	    corrosion depth&$\para_3$&0.25&1.00&\milli\meter \\
    \end{tabular}
    \caption{Parameters for the \textit{waveguide with a corrosion defect}.}
    \label{tab:CorroPara}
\end{table}

\begin{figure}
\centering
\subcaptionbox[Short Subcaption]{Objective function.\label{Fig:Corrosion_objective_function}}[0.45\textwidth ]
{\includegraphics[width=0.43\textwidth]{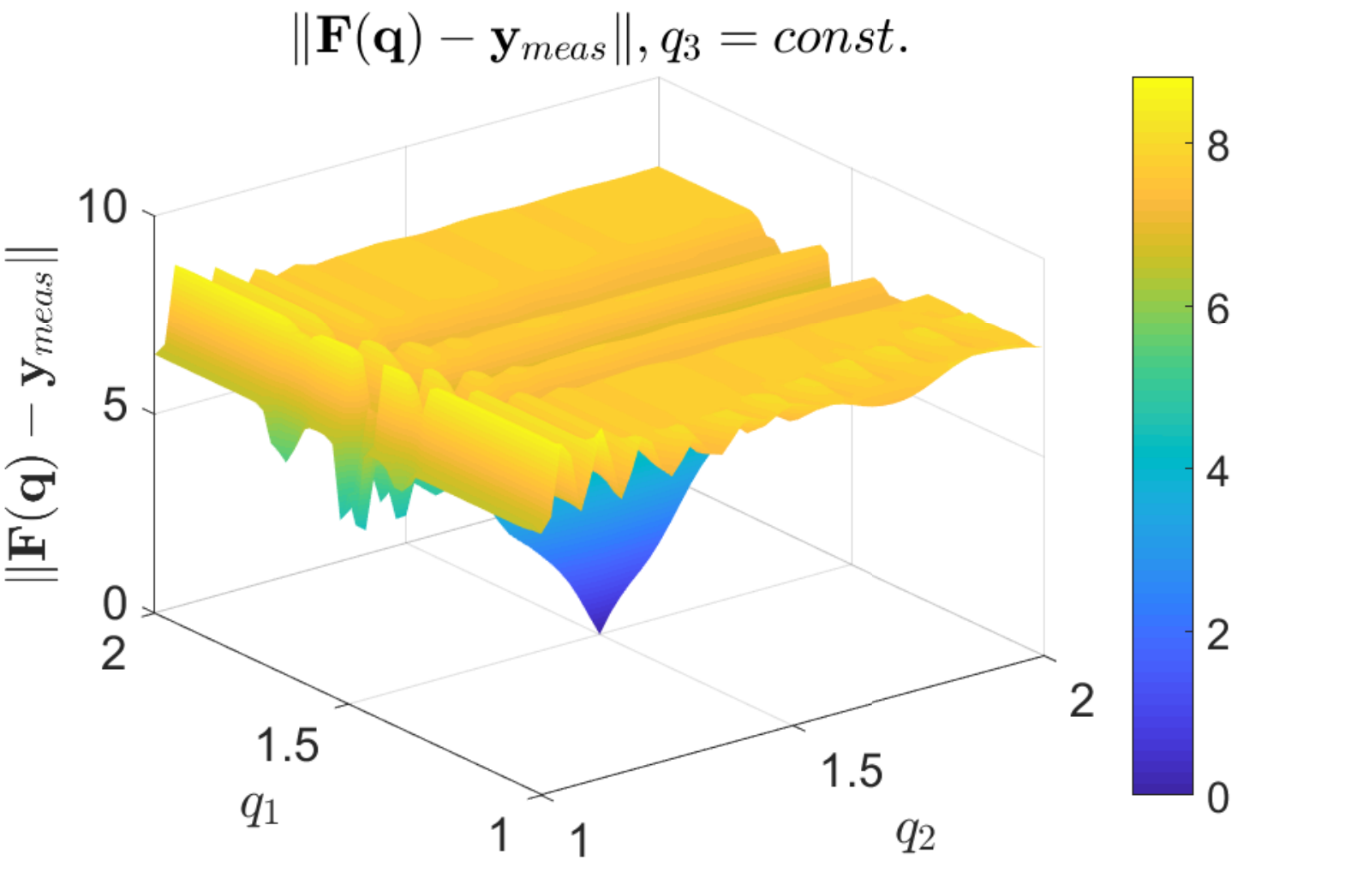} }
\hspace{0.05\textwidth} % seperation
\subcaptionbox[Short Subcaption]{Reconstruction error of parameters\\ $|\vpara_i^{\min}-\vpara_i^*|$ with noisy data. \label{fig:corrosion_plate_iterates_noise}}[0.45\textwidth]
{\includegraphics[width=0.43\textwidth]{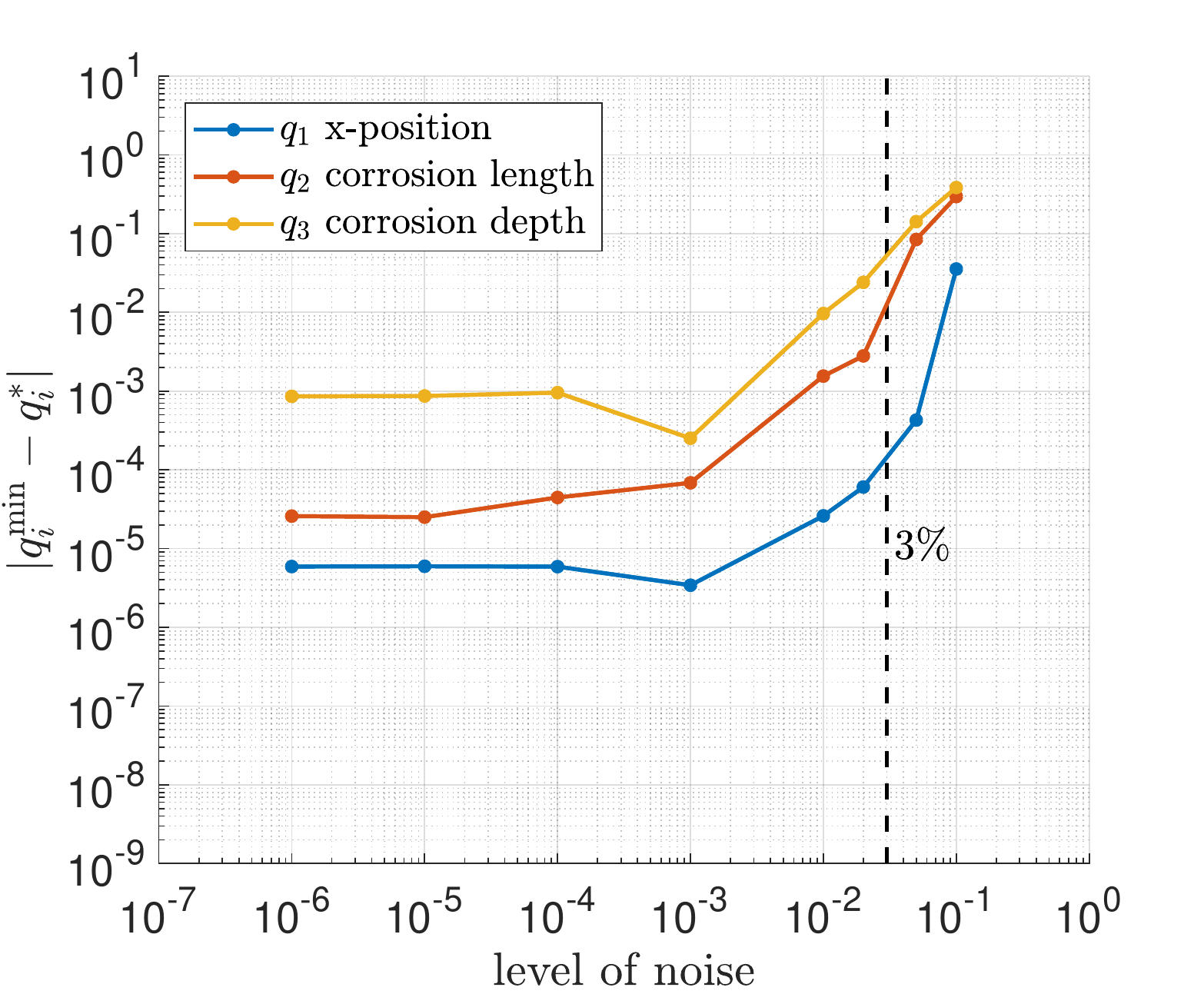} }
\caption[Short Caption]{
The objective function and reconstruction with noisy data for the \textit{waveguide with a corrosion defect}, where the target parameter $\vpara^*$ and the associated artificial measurement signal $\vy_{meas}$ correspond to the midpoint of the parameter space.}
\end{figure}

\begin{table}[htb]
\centering
\begin{tabular}{r|llllllllll}
$\para_1^*$& 1.59& 1.51& 1.47& 1.35& 1.55& 1.28& 1.52& 1.42& 1.34& 1.40 \\
$\para_2^*$& 1.43& 1.67& 1.28& 1.55& 1.46& 1.53& 1.44& 1.67& 1.29& 1.42 \\
$\para_3^*$& 1.48& 1.33& 1.44& 1.49& 1.57& 1.42& 1.52& 1.64& 1.54& 1.37 \\ \hline
$\|\vpara^{\min}-\vpara^*\|$& 9E-4& 9E-4& 9E-4& 8E-4& 9E-4& 8E-4& 9E-4& 9E-4& 9E-4& 9E-4
\end{tabular}
\caption{Error for reconstruction $\vpara^{\min}$ for differently shaped defects $\vpara^*\in [1,2]^3$ for the \textit{waveguide with a corrosion defect} with a noise level of $10^{-5}$.}
\label{tab:errors_corro}
\end{table}

The last example considers a steel plate with a 
simple model for a corrosion damage. The assumption is that parts of the steel are missing, and a tapering has formed in the waveguide.
Figure~\ref{fig:Corro} shows an overview of the model. Note that again the $x$-axis is interrupted. The traction $\ff_{S0}$ is the same as in first example (Equation~\eqref{eq:tractionS0}), which leads to an S0-mode excitation. The S0-mode shows greater sensitivity than the A0-mode to changes in thickness. Therefore, the S0-excitation is used.
The two re-entrant corners~\cite{bulling2019comparison} in sub-domains $\Omega_5$ and $\Omega_9$ are located at the scaling center to handle the possible singular stresses for a deep corrosion.
The model has $210$ degrees of freedom in total.

There are three parameters to be optimized. Once more, the first parameter defines the global position of the defect.
The second parameter defines the length of the corrosion by changing the size of the super element $\Omega_7$. 
The third parameter denotes the depth of the corrosion - see Figure~\ref{fig:Corro}.
The maximal and minimal values for the parameters are listed in Table~\ref{tab:CorroPara}.
Figure~\ref{Fig:Corrosion_objective_function} shows the objective function for a constant corrosion depth. The surface is more complicated due to the two main reflection points at the scaling centers, leading to mode conversion.
As mentioned earlier, we had to increase the number of random vectors drawn for the initial guess because of the more complicated shape of the objective function.
In Figure~\ref{fig:corrosion_plate_iterates_noise} the error for reconstruction with noise is given. 
Additionally, Figure~\ref{Fig:iterates_corro} shows the convergence behavior for the reconstruction error at the midpoint of the parameter space.
Despite the more complicated behavior of the objective function, the parameters show robustness against noise.
For all ten tests, the reconstruction leads to small errors as listed in Table~\ref{tab:errors_corro}.
	\section{Conclusion}\label{sec:Conclusion}

In the paper, we showed that a semi-analytical waveguide model given by SBFEM can be combined successfully with derivative-based optimization to reconstruct defects of different nature in a two-dimensional waveguide, even in the presence of noise. For this purpose, SBFEM is implemented in MATLAB, coupled with the AD tool ADiMat and embedded in an iterative solution procedure. The reconstruction can be performed on a standard desktop computer within a very moderate time such that our approach offers considerable runtime advantages compared to alternative methods. Three numerical tests were conducted to illustrate our approach: The first example covers a crack identification, the second one a delamination, and the third one a corrosion. In all three cases, the waveguides were made out of steel.
In the cross-sectional models, the reconstruction with a minimal amount of data is possible. Only the envelope of the $y$-displacement at a single point is sufficient if a certain geometry of a single defect is assumed. Due to the limited number of local minima, the envelope curve shows the most promising features for an optimization. In all tests carried out, the approach demonstrates robustness against noise and varying defect parameters.

Future work will be dedicated to the extension to the three-dimensional case and anisotropic material behavior. The future work will allow the use of the proposed approach for non-destructive testing and structural health monitoring in real-life applications like carbon-enforced structures.

\section*{Acknowledgments}
The authors gratefully acknowledge the German Research Foundation for funding (DFG project number 428590437).
	\printbibliography
	\appendix
	\section{Implementation Details}\label{sec:ID}

This appendix shows how the Lyapunov equation can be solved with the help of an eigenvalue decomposition.
The Lyapunov equation is
% A*X + X*A' + Q = 0
\begin{align}
	\mA \mX + \mX \mA^\ct + \mC = \vo. \label{eq:Lyapunov}
\end{align}
% A*V = V*D
% W'*A = D*W'
Define the right eigenvalue decomposition of $\mA$
\begin{align}
	\mV^{-1} \mA \mV & = \mD
\end{align}
with the right eigenvector matrix $\mV$ and the diagonal eigenvalue matrix $\mD$. Note that this leads to the left eigenvalue decomposition of $\mA^\ct$
\begin{align}
	\mV^{\ct} \mA^{\ct} \mV^{-\ct} & = \mD^{\ct}.
\end{align}
% CC = VA\(C*VB);
% X  = (VA*XX)/VB
After pre-multiplying with $\mV^{-1}$ and post-multiplying with $\mV^{-\ct}$, Equation~\eqref{eq:Lyapunov} leads to a new Lyapunov equation with a diagonal coefficient matrix
\begin{align}
	\mD \tilde \mX + \tilde \mX \mD^\ct + \tilde \mC = \vo
\end{align}
%XX = -CC./(repmat(dA,1,N1)+repmat(conj(dA).',N1,1));
with
\begin{align}
	\tilde \mC & = \mV^{-1} \mC \mV^{-\ct} && \text{and} & 
	\tilde \mX & = \mV^{-1} \mX \mV^{-\ct}. \label{eq:mX}
\end{align}
The matrix $\tilde \mX$ can be computed element-wise by
\begin{align}
	(\tilde \mX)_{ij} = \frac{- (\tilde \mC)_{ij}}{(\mD)_{ii} + (\mD)^\ct_{jj}}.
\end{align}
The solution matrix $\mX$ of the original Lyapunov equation can be derived by the inverse of Equation~\eqref{eq:mX}, i.e., $\mX = \mV^{\ct} \tilde \mX \mV$.

% Eigenvalue derivative
% \begin{align}
% 	- \mZ \mV & = \mV \mLbd \\
% 	- \mZ^\prime \mV - \mZ \mV^\prime & = \mV^\prime \mLbd + \mV \mLbd^\prime 
% \end{align}

% \begin{align}
% 	 \mV^\prime & = \mV \mC 
% \end{align}

% \begin{align}
% 	(\mLbd^\prime)_{ii} & = -(\mV^{-1} \mZ \mV)_{ii} \\
% 	(\mC)_{ij} & =
% 	\begin{cases}
% 		\frac{(\mV^1 \mZ)_{ij}}{\lambda_i - \lambda_j} &
% 		\lambda_i \neq \lambda_j\\
% 		0 & \lambda_i = \lambda_j
% 	\end{cases}
% \end{align}

	\newpage

\end{document}